\documentclass[reqno]{amsart}
\usepackage{amsmath, amsthm, amssymb, amstext, verbatim} 

\usepackage{hyperref,xcolor}%
\hypersetup{
 pdfborder={0 0 0},
 colorlinks=true,
 linkcolor=blue,
 citecolor=blue,
}

\usepackage[normalem]{ulem}

\usepackage[left=3cm,right=3cm,top=3cm,bottom=3cm]{geometry}

\usepackage[colorinlistoftodos]{todonotes}
\usepackage{enumitem}
\newtheorem{theorem}{Theorem}[section]
\newtheorem{proposition}[theorem]{Proposition}
\newtheorem{lemma}[theorem]{Lemma}
\newtheorem{e-proposition}[theorem]{Proposition}

\newtheorem{definition}[theorem]{Definition\rm}
\newtheorem{remark}[theorem]{Remark}

\def\d{\,{\rm d}}

\DeclareMathOperator*{\divergenz}{div}              %

\newcommand{\N}{\mathbb{N}}
\newcommand{\Z}{\mathbb{Z}}

\newcommand{\R}{\mathbb{R}}

\newcommand{\eps}{\varepsilon}

\newcommand{\Om}{\Omega}
\newcommand{\rand}{\partial\Omega}

\newcommand{\into}{\int_{\Omega}}

\newcommand{\bbu}{\boldsymbol{u}}
\newcommand{\bM}{{M}}
\newcommand{\bmu}{\boldsymbol{\mu}}

\numberwithin{theorem}{section}
\numberwithin{equation}{section}

\newcommand{\cc}{\mathfrak{C}}

\renewcommand{\d}[1]{\mathrm{d}#1}

\newcommand{\bfu}{\textbf{u}}

\newcommand{\bphi}{\boldsymbol{\phi}}

\newcommand{\bbw}{\boldsymbol{w}}

\newcommand{\A}{{\mathcal A}}

\title[Cross-diffusion equations with nonlocal Cahn-Hilliard terms]{Existence and local asymptotics for a system of cross-diffusion equations with nonlocal Cahn-Hilliard terms}

\author[E.\,Davoli]{Elisa Davoli}
\address[E.\,Davoli]{Tu Wien, Institute for Analysis and Scientific Computing, Wiedner Hauptstra\ss{}e 8-10, 1040 Wien, Austria}
\email{elisa.davoli@tuwien.ac.at}

\author[G.\,Marino]{Greta Marino}
\address[G.\,Marino]{Universit\"{a}t Augsburg, Institut f\"ur Mathematik, Universit\"{a}tsstra\ss e 12a, 86159 Augsburg, Germany}
\email{greta.marino@uni-a.de}

\author[J.-F.\,Pietschmann]{Jan-Frederik Pietschmann}
\address[J.-F.\,Pietschmann]{Universit\"{a}t Augsburg, Institut f\"ur Mathematik, Universit\"{a}tsstra\ss e 12a, 86159 Augsburg, Germany and Centre for Advanced Analytics and Predictive Sciences (CAAPS), University of Augsburg,
Universit\"{a}tsstr. 12a, 86159 Augsburg, Germany.}
\email{jan-f.pietschmann@uni-a.de}

\author[N.\,Zamponi]{Nicola Zamponi}
\address[N.\,Zamponi]{Universit\"{a}t Augsburg, Institut f\"ur Mathematik, Universit\"{a}tsstra\ss e 12a, 86159 Augsburg, Germany}
\email{nicola.zamponi@uni-a.de}



\begin{document}

\begin{abstract}

We study a nonlocal Cahn-Hilliard model for a multicomponent mixture with cross-diffusion effects and degenerate mobility. The nonlocality is described by means of a symmetric singular kernel.
We define a notion of weak solution adapted to possible degeneracies and prove, as our first main result, its  global-in-time existence. The proof relies on an application of the formal gradient flow structure of the system (to overcome the lack of a priori estimates),  combined with  an extension of the  boundedness-by-entropy  method, involving a careful analysis of an auxiliary variational problem. 
This allows to obtain solutions to an approximate, time-discrete system. Letting the time step size go to zero, we recover the desired nonlocal weak solution where, due to their low regularity, the Cahn-Hilliard terms require a special treatment.
Finally, we prove convergence of solutions for this class of nonlocal Cahn-Hilliard equations to their local counterparts.

\end{abstract}

\maketitle

\section{Introduction}

In this work we study a nonlocal Cahn-Hilliard model with degenerate mobility for a multicomponent mixture where cross-diffusion effects between the different species of the system are taken into account, and where the species do separate from each other. The motivation for considering such a model stems from multiphase systems modelling isothermal phase separation of  miscible entities occupying an isolated region $\Om \subset \R^d$ ($d= 1,2,3$), cf. \cite{Elliott-Luckhaus} and the references therein. 
We consider a cross-diffusion system that models the interactions between $n+1$ species, $n \in \N \setminus \{0\}$. %
More precisely, let {$\Omega$ be the $d$-dimensional flat torus of dimension $d = 2,3$} (which corresponds to imposing periodic boundary conditions), and let $T>0$ be some final time.
We assume that the $n+1$ species in the mixture occupy the spatial domain $\Omega$ and for all $i= 0, \dots, n$, we denote by $u_i(t,x)$ the volume fraction of the $i^{th}$ species at time $t\in [0,T]$ and point $x\in \Omega$ and set $\bbu:= (u_0, u_1, \dots, u_n)$. Given a small parameter $\eps>0$ which accounts for the  radius of the  nonlocal interactions, we denote by $\bbu_\eps= (u_{\eps,0}, \dots, u_{\eps, n})$ the nonlocal counterpart  to  $\bbu$. 

We are then interested in the existence of weak solutions to the following system of cross-diffusion equations with nonlocal interactions
        \begin{equation}
        \label{intro1}
    \partial_t \bbu_{\eps}= \divergenz \left(\bM(\bbu_{\eps}) \nabla \bmu_{\eps}\right), 
    \end{equation}
such that 
    \[
      0\leq u_{\eps,i}(t,x) \leq 1 \text{ for  every  } i=0, \dots, n   \text{ and } \sum_{i=0}^n u_{\eps,i}(t,x) = 1\text{ 
for a.e. }(t,x)\in [0,T]\times \Omega.
    \]
Here, for all $\bbu_\eps \in \mathbb{R}_+^{n+1}$, $ M(\bbu_\eps)\in \mathbb{R}^{(n+1) \times (n+1)}$ is a degenerate mobility matrix whose precise expression is given in Section~\ref{sect2}, while $\bmu_\eps$ is the chemical potential, defined as 
\begin{equation*} %
\bmu_\eps = D_{\bbu_\eps} E_{ NL}(\bbu_\eps),
\end{equation*}
for  $E_{\text{NL}}$
given by
    \begin{equation}
        \label{energy-NL}
        \begin{aligned}
    E_{\text{NL}}({\bbu}_{\eps}):&= \sum_{i=0}^n \into u_{\eps,i} \ln u_{\eps, i}- u_{\eps, i}+ 1  {-\frac{1}{2}\sum_{i,j=1}a_{ij}u_{\eps, i}u_{\eps, i}} \; \d x\\
    &+ \frac{1}{4} \sum_{i,j=0}^n \into \into c_{ij} K_{\eps}(x, y) (u_{\eps, i}(x)- u_{\eps, i}(y)) (u_{\eps, j}(x)- u_{\eps, j}(y)) \; \d x \d y,
    \end{aligned}
    \end{equation}
where  $K_{\eps}$ are $L^1$ and symmetric convolution kernels and $ {\cc}= (c_{ij})_{i,j=0,\ldots,n} \in \R^{(n+1) \times (n+1)}$

is a positive semidefinite matrix, see Section \ref{sect2} for their definitions. The matrix $A= (a_{ij})_{ij} \in \R^{(n+1) \times (n+1)}$ is symmetric and positive semi-definite. Therefore, the term $-\frac{1}{2}\bfu\cdot A\bfu$ is concave. Note that the logarithmic terms in this energy functional account for diffusion, {the second term induces separation of the phases and the } integral terms are responsible for nonlocal {interaction}. In the case when $K_\eps \equiv 0$ {and $A=0$}, system (\ref{intro1}) boils down to a multi-species degenerate cross-diffusion system with size exclusion that was studied for example in~\cite{BDPS2010,BE2018freeboundary,Berendsen2019,cances2020convergent}. 
The local counterpart  to  \eqref{intro1}, namely system
        \begin{equation}
        \label{local-syst}
    \partial_t \bbu= \divergenz \left(\bM(\bbu) \nabla \bmu\right), 
    \end{equation}
{is driven by }the local energy {functional}  {$E_{\text{L}}$}
given by
\begin{equation}
\label{local-energy}
E_{\rm L}({\bbu}):= \sum_{i=0}^n \into 
(u_i \ln u_i- u_i+ 1)\d x 
{-\frac{1}{2}\sum_{i,j=1}\into a_{ij}u_iu_j\d x} + \frac{1}{2} \sum_{i,j=0}^n \into c_{ij} \nabla u_i \cdot \nabla u_j  \; \d x.
\end{equation}
Indeed, as $\eps \to 0$, it  formally  holds that
    \[
    \frac{1}{4} \sum_{i,j=0}^n \into \into c_{ij} K_{\eps}(x, y) (u_{i}(x)- u_{i}(y)) (u_{j}(x)- u_{j}(y)) \; \d x \d y \to \frac{1}{2} \sum_{i,j=0}^n \into c_{ij} \nabla u_i \cdot \nabla u_j  \; \d x.
    \]

The main goals of our paper are to show existence of solution to the non-local system \eqref{intro1} as well as convergence of the non-local to the local equations. Before we proceed, let us put our work into perspective with respect to previous results.

\subsubsection*{Cross-diffusion systems with size exclusion} 

Systems of partial differential equations with cross-diffusion have gained a lot of interest in recent years \cite{Kfner1996InvariantRF,Chen2004,Chen2006,Lepoutre2012,juengel} 
and appear in many applications, for instance the modeling of population dynamics of multiple species~\cite{Burger2016} or cell sorting or chemotaxis-like applications~\cite{Painter2002,Painter2009}. One major difficulty in the analysis of such strongly coupled systems is the lack of a priori estimates guaranteeing enough compactness of the maps $u_{\eps,i}$,  an issue that can be overcome by the \emph{boundedness-by-entropy method} if the system has the so-called \emph{formal gradient flow} structure, see  \cite{BDPS2010, juengel} for a more detailed presentation.

\subsubsection*{Cross-diffusion systems with Cahn-Hilliard contributions} 
These kind of systems have been recently studied in \cite{EMP}, {see also \cite{juengel}  for a different choice of mobility}, and describe the evolution of a multicomponent mixture where cross-diffusion effects between the different species are taken into account, and (in this particular case) where only one species does separate from the others. This is motivated by multiphase systems where miscible entities may coexist in one single phase, see \cite{klinkert2015comprehension,Wenisch2016} for examples. Within this phase, cross-diffusion between the different species is taken into account in order to correctly account for finite size effects that may occur at high concentrations.

\subsubsection*{Nonlocal Cahn-Hilliard equations} A nonlocal model for phase separation was originally proposed in \cite{GL}. Ever since the works \cite{MRST,DRST}, the literature on nonlocal-to-local convergence of Cahn-Hilliard models has bloomed. We refer to \cite{DGL,AH,ES} for an overview of the most recent developments. For well-posedness, regularity, and strict separation properties of multi-component Cahn-Hilliard systems with singular potentials we refer to the recent work \cite{Gras}.

\medskip

\subsubsection*{Nonlocal cross-diffusion equations}

Nonlocal interaction equations are continuum models for large systems of particles where every single particle can interact not only with its immediate neighbors but also with particles far away. These equations have a wide range of applications. In biology they are used to model the collective behavior of a large number of individuals \cite{ME}, in physics they are used in models describing the evolution of vortex densities in superconductors \cite{We}, and they appear in simplified inelastic interaction models for granular media \cite{BCP}. In their simplest form, nonlocal interaction equations can be written as
	\begin{equation}
	\label{nonlocal}
	\frac{\partial \sigma}{\partial t}+ \operatorname{div} (\sigma r)= 0, \quad r= -\nabla W \ast \sigma,
	\end{equation}
where $\sigma(t, x)$ is the probability or a mass density of particles at time $t$ and at location $x \in \mathbb R^N$, $W\colon \mathbb R^N \to \mathbb R$ is the interaction potential and $r(t, x)$ is the velocity of the particles. 
It is by {now well understood, see for example \cite{ambrosioGradientFlowsMetric2008}, }
that equation \eqref{nonlocal} is a gradient flow of the interaction energy
	\[
	E(\sigma)= \frac{1}{2} \iint_{\mathbb R^N \times \mathbb R^N} W(x- y) \; d\sigma(x) d\sigma(y)
	\]
with respect to the {$2$-}Wasserstein distance. Thus, stable steady states of \eqref{nonlocal} are expected to be local minimizers of the interaction energy. 
Much  less  is  known  however  in  the  case  of  multi-species  systems,  which  received  most attention  only  recently.   With  different  species, the  modeling  leads  to  nonlinear  degenerate  cross-diffusion  systems  for  the  densities  of  all species,  again  with  some  nonlocal  terms.  First rigorous studies of stationary problems   show  interesting  phase  separation  phenomena,  whose  dynamics seems rather unexplored so far \cite{BSW12, Cic16}.  In \cite{Berendsen2017nonlocal} the authors study a nonlocal cross-diffusion model for two species, which can be derived from a lattice-based microscopic model with size exclusion. In  order  to  see  the  inherent  phase  separation  in their model, they  rewrite  it  as  a  system  of nonlocal Cahn-Hilliard equations by means of a nonlocal Laplacian  defined as a negative semidefinite operator. 
 For the complementary setting in which the diffusion is driven by a local operator but further nonlocalities arising, e.g., in the mean-field limit of stochastic PDEs in the neurosciences, we refer to the two recent contributions \cite{JPZ,JPZ2}. 

\subsubsection*{Degenerate cross-diffusion systems with nonlocal terms} Concerning nonlinear degenerate systems with non-local terms, very little is known on the dynamics, which seems rather unexplored so far \cite{BSW12, Cic16}. In \cite{Berendsen2017nonlocal} the authors study a nonlocal cross-diffusion model for two species, which can be derived from a lattice-based microscopic model with size exclusion. In order to see the inherent phase separation in their model, they rewrite it as a system of nonlocal Cahn-Hilliard equations by means of a nonlocal Laplacian defined as a negative semidefinite operator.

\subsection*{Contribution and structure of the paper}
In this article we combine ideas mainly coming from \cite{DRST} and \cite{EMP} and prove the existence of global weak solutions to system \eqref{intro1} with energy \eqref{energy-NL} and supplemented with appropriate initial- and boundary conditions and then show that this solution, in the limit $\eps \to 0$, converges to a solution to the local counterpart of \eqref{intro1}.\medskip\\
 The novelty of our work  is threefold. 
 
	\begin{itemize}
 
		\item[(a)] This is, to the best of our knowledge, the first attempt to combine degenerate cross-diffusion system with non-local Cahn-Hilliard contributions. 
  
		\item[(b)] We are able to treat an energy that  involves non-local Cahn-Hilliard terms acting on all species.  Note that this requires  an appropriate definition of weak solutions and a careful analysis when performing the limit of an approximate time-discrete non-local system, {as the non-local terms yield less regularity in the a priori estimates than their local counterparts.}

  		\item[(c)] We apply the generalized boundedness-by-entropy method exploited in \cite{EMP} to this non-local context.
	\end{itemize}

\medskip

This manuscript is organized as follows. In Section~\ref{sect2}, we  provide the setting of the problem, introduce our notion of weak solutions and state the main existence theorems. The proof in the non-local case is based on the introduction of a regularized time discrete approximate problem, depending on a positive time step $\tau$, which is presented in Section~\ref{sect3}. We derive a priori estimates and prove the existence of time-discrete iterates via a Schauder fixed point argument. Then we recover some regularity properties and exploit them to pass to the limit as the time step $\tau \to 0$ and obtain a solution to the nonlocal system \eqref{intro1}. Finally, in Section \ref{sect4}, we show that the a priori estimates previously obtained actually do not depend on the non-local parameter $\eps$, which gives us enough regularity to pass to the limit as $\eps \to 0$ and recover a solution to the local system \eqref{local-syst}.

\section{ Setting of the problem and main results}
\label{sect2}

In this section we  specify our setting, provide the notions of non-local and local weak solutions, and state the main results. 

\subsection{Formulation of the model}

Let $\Omega{=\mathbb{T}^d}$ be the $d$-dimensional flat torus, $d=  2, 3$, (this corresponds to imposing periodic boundary conditions)  and let $u_i= u_i(t,x)$, $i= 0, 1, \dots, n$, be the volume fractions of  different atomic species that occupy the volume of $\Omega$ at the time $t>0$.  In what follows, we will always implicitly identify all the maps with their restrictions to $\Omega$, automatically resulting into periodic functions. Additionally, with a slight abuse of notation, we will identify the unit torus with the unit cube ${\left[-\frac{1}{2}, \frac{1}{2}\right]^d}$ in $\R^d$. Furthermore, for every $\eps>0 $ we consider a family of kernels $K_\eps : \R^d\times\R^d\to\R_+$ satisfying $K_\eps(x,y)=\tilde K_\eps(x-y)$ $\forall x,y\in\R^d$, for a function $\tilde K_\eps : \R^d\to\R_+$ which is $1$-periodical with respect to each coordinate:
$$
\tilde K_\eps(z+e^{i}) = \tilde K_\eps(z)\quad\forall z\in \R^d,\mbox{ with }e^i\in\R^d:~~(e^i)_j = \delta_{ij},~~i,j=1,\ldots,d,
$$
and admits the following representation on the unit cube $\left[-\frac{1}{2}, \frac{1}{2}\right]^d$:
\begin{equation}
\label{kernel}
\tilde K_{\eps}(z):= \frac{\rho_{\eps}(|z|)}{|z|^2},\quad\forall z\in \left[-\frac{1}{2}, \frac{1}{2}\right]^d,
\end{equation}
for a suitable family $(\rho_{\eps})_{\eps}$ of $C^{\infty}$-maps satisfying the following properties (see \cite{Ponce,Ponce2}): 
{$x\in\R^d\mapsto\rho_\eps(|x|)\in\R$ is $C^\infty$,}
$\rho_{\eps}:\R\to[0,+\infty)$ is defined as $\rho_{\eps}(r):=\frac{1}{\eps^d}\rho(r/\eps)$, with 
$\rho\in C^{\infty}_c({[0,\frac{1}{2})})$ 
such that
$$
{\alpha_d\equiv 
\int_{\R^d}\frac{\rho(|x|)}{|x|^2}\,dx = \omega_d}
\int_0^{+\infty}\rho(r)r^{d-3}\,dr<+\infty,\quad 
{\int_{\R^d}\rho(|x|)\,dx = \omega_d}
\int_0^{+\infty}\rho(r)r^{d-1}\,dr=2d,
$$
{where $\omega_d$ is the $(d-1)$-dimensional Hausdorff measure of the sphere $\mathbb{S}^{d-1}$, namely $\omega_2 = 2\pi$, $\omega_3 = 4\pi$.}
Note that under these assumptions $K_{\eps}\in L^1(\Omega\times \Omega)$ and $K_{\eps}(x,y)=K_{\eps}(y,x)$ for almost every $x,y\in \Omega$.

Set $\bbu_{\eps}= \bbu_{\eps}(t,x):= (u_{\eps,0}, u_{\eps,1}, \dots, u_{\eps,n})$. We are interested in the existence of weak solutions to the following system 
	\begin{equation}
	\label{system}
 	\begin{split}
	\partial_t \bbu_{\eps}&= \divergenz \left(\bM(\bbu_{\eps}) \nabla \bmu_{\eps}\right), 
	\end{split}
	\end{equation}
 {where the above expression is to be understood in the sense that $\partial_t u_{\eps,i} = \divergenz( M(\bbu_{\eps})\nabla \mu_{\eps})_i$ for $i=0,\ldots, n$, and }such that
	\begin{equation}
 \label{constraint}
	0 \le u_{\eps, i}(t,x) \le 1 \text{ and } \sum_{i=0}^n u_{\eps, i}(t,x)= 1\quad \text{for all } i= 0, \dots, n  \text{ and a.e. } (t,x) \in [0,T] \times \Om.
	\end{equation}
Note that $\bM(\cdot){\in \R^{(n+1)\times (n+1)}}$ in \eqref{system} is a degenerate mobility matrix with components
	\begin{equation}
	\label{mobility}
	\begin{aligned}
	M_{ij}({\bbu}_{\eps})&:= -L_{ij} u_{\eps, i} u_{\eps, j} \quad && \forall \, i \ne j= 0, \dots, n, \\
	M_{ii}({\bbu}_{\eps})&:= \sum_{0 \le j \ne i \le n} L_{ij} u_{\eps, i} u_{\eps, j} && \forall \, i= 0, \dots, n,
	\end{aligned}
	\end{equation}
with $L_{ij}>0 $ such that $L_{ij}= L_{ji}$, for all $i \ne j= 0, \dots, n$, while $\bmu_{\eps}$ is the chemical potential, defined formally as
	\[
	\bmu_\eps:= D_{\bbu_\eps} E_{ NL}(\bbu_\eps),
	\]
 where  {$D_{\bbu_\eps} E_{\text{NL}}(\bbu_\eps)$ denotes the variational derivative of the functional }
	\begin{equation}
	\label{E}
	 E_{ NL}({\bbu}_{\eps}):= \into F(\bbu_{\eps}) \; \d x+ \frac{1}{4} \sum_{i,j=0}^n \into \into c_{ij} K_{\eps}(x, y) (u_{\eps, i}(x)- u_{\eps, i}(y)) (u_{\eps, j}(x)- u_{\eps, j}(y)) \; \d x \d y,
	\end{equation}
 with 
	\[
	F(\bbu_{\eps})= \sum_{i=0}^n (u_{\eps, i} \ln u_{\eps, i}- u_{\eps, i}+ 1) {-\frac{1}{2}\sum_{i,j=0}^n a_{ij}u_{\eps, i}u_{\eps, j}}.
	\]
This in particular gives 
	\begin{equation}
	\label{mu-i}
	\begin{split}
	\mu_{\eps, i}= D_{u_{\eps, i}}{ E_{NL}({\bbu}_{\eps})}&= \ln u_{\eps, i}{-\sum_{j=0}^n a_{ij}u_{\eps, j}}+ \sum_{j=0}^n c_{ij} \big((K_{\eps} \ast 1) u_{\eps, j}- K_{\eps} \ast u_{\eps, j}\big) \\
	&=: \ln u_{\eps, i}{-\sum_{j=0}^n a_{ij}u_{\eps, j}}+ \sum_{j=0}^n c_{ij} B_{\eps}(u_{\eps, j}),
	\end{split}
	\end{equation}
for all $i=0, \dots, n$, so that  $\bmu_{\eps}= (\mu_{\eps, 0}, \dots, \mu_{\eps, n})$,  while the properties of the operator 
\begin{align}\label{Beps}
{B_{\eps}(v) := (K_{\eps} \ast 1) v- K_{\eps} \ast v}
\end{align}
are as in \cite[Section 2.2]{DRST}. Note that  we set 
    \[
    (K_\eps \ast 1)(x)= \into K_\eps(x,y) \; \d y \quad \text{as well as} \quad (K_\eps \ast v)(x)= \into K_\eps(x,y) v(y) \; \d y.
    \]
{We point out that via Young's inequality for convolutions it follows
\begin{align}
\label{Beps.esteps}
\|B_\eps(u)\|_{L^p(\Omega)}\leq 2\alpha_d\eps^{-2}\|u\|_{L^p(\Omega)},\quad\forall u\in L^p(\Omega),\quad\forall p\in [1,\infty].
\end{align}
}
Under the previous assumptions on the interaction kernels, {an adaptation of} the argument  in \cite[Proof of Theorem 2.2]{DRST} shows  that the following result holds true.
\begin{lemma}\label{lem:laplace}
Let $(v_\eps)_{\eps}\subset H^1(0,T;(H^1)(\Omega)')\cap L^2(0,T;H^1(\Omega))$ be such that
$$\|v_\eps\|_{H^1(0,T;(H^1)(\Omega)')\cap L^2(0,T;H^1(\Omega))}+\|B_{\eps}(v_{\eps})\|_{L^2(0,T;{L^2(\Omega)})}\leq C.$$

Then, up to subsequences, 
\begin{align*}
    & v_{\eps}\to v \quad \text{weakly in }H^1(0,T;(H^1)(\Omega)')\cap L^2(0,T;H^1(\Omega)),\\
    &B_{\eps}(v_{\eps})\to -\Delta v\quad \text{weakly in }L^2(0,T;{L^2(\Omega)}).
\end{align*}
\end{lemma}
{
Let us preliminarily observe that, given any smooth functions $u_1, u_2 : \Omega\to\R$, it holds
$$
B_\eps(\nabla u_1) (x) = \int_{\Omega}\frac{\rho_\eps(|x-y|)}{|x-y|^2}(\nabla u_1(x)-\nabla u_1(y))\d y,
$$
hence it follows via symmetry
\begin{align}\label{B.qf}
\int_{\Omega}\nabla u_2\cdot B_\eps(\nabla u_1)\d x = 
\frac{1}{2}\int_{\Omega\times\Omega}K_\eps(x,y)
(\nabla u_1(x)-\nabla u_1(y))\cdot 
(\nabla u_2(x)-\nabla u_2(y))\d x\d y .
\end{align}
We are now going to show the following inequality, which will be used in the nonlocal to local limit:
\begin{lemma}\label{lemma.new.B.eps}
There exists a constant $C>0$ independent of $\eps>0$ such that
\begin{align}
\label{new.1}
\|B_\eps(v)\|_{L^2(\Omega)}^2\leq C\int_{\Omega\times\Omega}K_\eps(x,y)|\nabla v(x)-\nabla v(y)|^2\d x\d y, \qquad\forall v\in H^1(\Omega).
\end{align}
\end{lemma}
\begin{proof}
We can assume $v\in C^\infty(\Omega) = C^\infty(\mathbb{T}^d)$ and then argue by density.
From \eqref{B.qf} we deduce that \eqref{new.1} is equivalent to
\begin{align}
\label{new.2}
\|B_\eps(v)\|_{L^2(\Omega)}^2\leq C\left(
\nabla v, B_\eps(\nabla v)
\right)_{L^2(\Omega)},\qquad\forall v\in C^\infty(\Omega).
\end{align}
Let $\mathcal{F}$ be the discrete Fourier transform. From the definition of $B_\eps$ it follows
\begin{align*}
\mathcal{F}(B_\eps v)(k) = \hat{v}(k)\int_{\mathbb{T}^d}\frac{\rho_\eps(|z|)}{|z|^2}\d z - 
\hat{v}(k)\mathcal{F}(|\cdot|^{-2}\rho_\eps(|\cdot|)(k) = 
b_\eps(k)\hat{v}(k),\qquad\forall k\in\Z^d,
\end{align*}
with
\begin{align*}
b_\eps(k) = \int_{\mathbb{T}^d}\frac{\rho_\eps(|z|)}{|z|^2}(1-e^{2\pi ik\cdot z})\d z = \int_{\mathbb{T}^d}\rho_\eps(|z|)\frac{1-\cos(2\pi k\cdot z)}{|z|^2}\d z,\qquad\forall k\in\Z^d.
\end{align*}
Thus the Fourier transform of the operator $B_\eps$ is simply the multiplication operator times the real, nonnegative, even function $b_\eps$. Therefore, via Parseval's identity inequality \eqref{new.2} is equivalent to
\begin{align}
\label{new.3}
\sum_{k\in\Z^d} b_\eps(k)^2 |\hat{v}(k)|^2\leq 
C\sum_{k\in\Z^d} |k|^2 b_\eps(k)|\hat{v}(k)|^2,
\qquad\forall v\in C^\infty(\Omega).
\end{align}
A sufficient condition for \eqref{new.3} to hold for every Schwartz function $v$ (and a fortiori for every $v\in C^\infty(\Omega)$) is that the following inequality
\begin{align}\label{new.4}
b_\eps(k)\leq C|k|^2,\qquad\forall k\in\Z^d,
\end{align}
be verified. This is indeed the case, since
\begin{align*}
b_\eps(k) =& |k|^2 \int_{\mathbb{T}^d}\rho_\eps(|z|)\frac{1-\cos(2\pi k\cdot z)}{|k|^2 |z|^2}\d z\leq 
|k|^2 \int_{\mathbb{T}^d}\rho_\eps(|z|)\frac{1-\cos(2\pi k\cdot z)}{(k\cdot z)^2}\d z
\leq C|k|^2 \int_{\mathbb{T}^d}\rho_\eps(|z|)\d z\\ 
=& C|k|^2
\end{align*}
thanks to the normalization condition $\int_{\mathbb{T}^d}\rho_\eps(|z|)\d z=1$ and the property $\sup_{s\neq 0}\frac{1-\cos(2\pi s)}{s^2}<\infty$.
This finishes the proof of the Lemma.
\end{proof}
}

We will also make use of the following inequality:
\begin{align}\nonumber
&\forall\delta>0\quad \exists\eps_\delta>0,~~C_\delta>0:\\ 
\label{eq:h1_control_nonlocal}
&\|u\|_{H^1(\Omega)}^2 \le \delta \int_\Omega\int_\Omega K_\eps(x,y) 
|\nabla u(x)-\nabla u(y)|^2 \d x \d y + C_\delta \|u\|_{L^2(\Omega)}^2,\\ \nonumber
&\forall\eps<\eps_\delta,\quad 
\forall u\in H^1(\Omega):~~\int_\Omega\int_\Omega K_\eps(x,y) |\nabla u(x)-\nabla u(y)|^2 \d x \d y<\infty.
\end{align}
which can be easily proved by contradiction in a similar way as \cite[Lemma 4]{DRST} by employing \cite[Theorem 1.2]{Ponce2}.

Finally, the matrix ${\cc}= (c_{ij})_{i,j=0, \dots, n} $ 
is symmetric and there exists a constant $\eta_C>0$ such that 
    \begin{align}
    \label{lb.C}
    \xi^T {\cc} \xi \ge \eta_C \|\xi\|^2, \quad \text{ for all } \xi \in \R^{n+1}
    \end{align}

Then, system \eqref{system} can be written in the  scalar form
	\begin{align*}
	\partial_t u_{\eps, i}&= \divergenz \big((\bM({\bbu_\eps}) \nabla \bmu_{\eps})_i\big)
    = \divergenz \Big(\sum_{0 \le j \ne i \le n} L_{ij} u_{\eps, i} u_{\eps_j} \nabla (\mu_{\eps, i}- \mu_{\eps, j})\Big) \nonumber \\
    &= \divergenz \Big(\sum_{0 \le j \ne i \le n} L_{ij} \Big[u_{\eps, j} \nabla u_{\eps, i}- u_{\eps, i} \nabla u_{\eps, j}- u_{\eps, i} u_{\eps, j} \nabla (z_i(\bbu_\eps)-z_j(\bbu_\eps)) \\
    &\qquad\qquad\qquad\qquad\quad 
    + u_{\eps, i} u_{\eps, j} \nabla( q_{\eps,i}(\bbu_\eps)-  q_{\eps,j}(\bbu_\eps))\Big]\Big), \nonumber 
	\end{align*}
for all $i= 0, \dots, n$, where
    \begin{equation}
        \label{q}
    q_{\eps,i}(\bbu_\eps):= \sum_{k=0}^n c_{ik} B_{\eps}(u_{\eps, k}) \quad \forall \, i= 0, \dots, n.
    \end{equation}
    {and 
    \begin{equation}
        \label{z}
    z_i(\bbu_\eps):= \sum_{k=0}^n a_{ik} u_{\eps, k} \quad \forall \, i= 0, \dots, n.
    \end{equation}
    }
 For all $i \ne j= 0, \dots,n$ we {denote the part of the flux caused by the non-local contributions as}
    \[
    {J}_{\eps,ij}:= u_{\eps, i} u_{\eps,j} \nabla(q_{\eps,i}(\bbu_\eps)- q_{\eps,j}(\bbu_\eps)),
    \]
and set ${J}_\eps:= ({J}_{\eps,ij})_{ij}$. Finally, let $\bbu_\eps^0= (u_{\eps,0}^0, \dots, u_{\eps,n}^0)$ be an initial condition such that
	\begin{equation}\label{in-cond}
	u_{\eps,i}^0(x) \ge 0 \quad \text{ for all } i= 0, \dots, n, \quad \sum_{i=0}^n u_{\eps,i}^0(x)=1,  \quad \text{a.e. in } \Om.
	\end{equation}
For every $\eps>0$, we then  look for a solution $(\bbu_\eps, {J}_\eps)$ to the following problem
	\begin{equation}
	\label{scalar-nonlocal}
	\begin{aligned}
	& \partial_t u_{\eps, i}= \divergenz \Big(\sum_{0 \le j \ne i \le n} L_{ij} \Big[u_{\eps, j} \nabla u_{\eps, i}- u_{\eps, i} \nabla u_{\eps, j}- u_{\eps, i} u_{\eps, j} \nabla(z_i(\bbu_\eps)-z_j(\bbu_\eps)) + {J}_{\eps,ij}\Big]\Big) && \text{in } \Om \times (0,T), \\
	& \sum_{i=0}^n u_{\eps, i}= 1 && \text{in } \Om \times (0,T), \\
 & {J}_{\eps,ij}= u_{\eps, i} u_{\eps,j} \nabla(q_{\eps,i}(\bbu_\eps)- q_{\eps,j}(\bbu_\eps)) && \text{in } \Om \times (0,T), \\
	& u_{\eps, i}(\cdot \,, 0)= u_{\eps, i}^0 && \text{in } \Om,
	\end{aligned}
	\end{equation}
for all $i=0, \dots, n$. 
Similarly, {for the local system}, if $\bbu^0$ is an initial condition such that
    \[
    	u_{i}^0(x) \ge 0 \quad \text{ for all } i= 0, \dots, n, \quad \sum_{i=0}^n u_{i}^0(x)=1\text{ a.e. in } \Om,
    \]
by $\bbu$ we indicate a solution to the following {local} problem
	\begin{equation}
	\label{scalar-local}
	\begin{aligned}
	& \partial_t u_{i}= \divergenz \Big(\sum_{0 \le j \ne i \le n} L_{ij} \Big[u_{j} \nabla u_{ i}- u_{ i} \nabla u_{j}{- u_{i} u_{j} \nabla (z_i(\bbu)-z_j(\bbu))}- u_{ i} u_{ j} \nabla\sum_{k=0}^n (c_{ik}-c_{jk})\Delta u_{k}\Big]\Big) \quad && \text{in } \Om \times (0,T), \\
	& \sum_{i=0}^n u_{ i}= 1 && \text{in } \Om \times (0,T), \\
	& u_{ i}(\cdot \,, 0)= u_{ i}^0 && \text{in } \Om.
	\end{aligned}
	\end{equation}

\subsection{Notions of weak solution and statement of main results}

In this subsection we state the main results of the paper. The idea is to pursue the following goals:

\begin{enumerate}

\item To show that system \eqref{scalar-nonlocal} admits at least a weak solution ${\bbu}_{\eps}$ in a suitable sense that will be specified below;

\item To show that, since $(K_{\eps} \ast 1)v- (K_{\eps} \ast v) \to -\Delta v $ as $\eps \to 0^+$, solutions to \eqref{scalar-nonlocal} converge to solutions to \eqref{scalar-local}.

\end{enumerate}

\begin{definition}[Solution to the nonlocal system]
\label{def-nonlocal}

Let $\eps>0$ and $T>0$ be fixed and let $\bbu_\eps^0 \in H^1(\Om)$ be an  initial condition satisfying \eqref{in-cond}. We say that $(\bbu_\eps, {J}_\eps)= \big((u_{\eps,i})_{i \in \{0, \dots, n\}}, ({J}_{\eps, ij})_{ij \in \{0, \dots, n\}} \big)$ is a solution to the nonlocal system \eqref{scalar-nonlocal} if

	\begin{enumerate}
	
	\item[$(1_{\text{NL}})$] $0 \le u_{\eps, i} \le 1 $ for all $i=0, \dots, n$ and $\displaystyle \sum_{i=0}^n u_{\eps, i}= 1$ a.e. in $\Om \times (0,T)$; 
	
	\item[$(2_{\text{NL}})$] $u_{\eps, i} \in L^2(0,T; H^1(\Om))$ and $\partial_t u_{\eps, i} \in L^2(0,T; H^1(\Om)^\prime)$ for all $i=0, \dots, n$; 
	
	\item[$(3_{\text{NL}})$] 
    {$u_{\eps, i}(\cdot \,, t)\to u_{\eps, i}^0 $
    strongly in $L^2(\Omega)$ as $t\to 0$,}
    for all $i= 0, \dots, n$;

    \item[$(4_{\text{NL}})$] $ J_{\eps, ij} \in L^2((0,T) \times \Om)^d $;

    \item[$(5_{\text{NL}})$] ${J}_{\eps, ij}= u_{\eps,i} u_{\eps,j} \nabla(q_{\eps,i}(\bbu_\eps)- q_{\eps,j}(\bbu_\eps))$ in the following weak sense
        \[
        \int_0^T \into {J}_{\eps,ij} \cdot {\eta} \; \d x \d t= - \int_0^T \into \divergenz(u_{\eps,i} u_{\eps,j} {\eta}) (q_{\eps,i}(\bbu_\eps)- q_{\eps,j}(\bbu_\eps)) \; \d x dt
        \]
    for every ${\eta} \in L^2(0,T; H^1(\Om)^d) \cap L^\infty ((0,T) \times \Om; \R^d) $, with ${\eta} \cdot {n}= 0$ on $\rand$, {where $q_{\eps,i}$ is defined in \eqref{q}};
	
	\item[$(6_{\text{NL}})$] for every $i= 0, \dots, n $ and every $\varphi_i \in L^2(0,T; H^1(\Om))$ there holds
    \begin{equation}
    \label{eq:nonlocal}
    		\begin{split}
		& \int_0^T \langle \partial_t u_{\eps, i}, \varphi_i \rangle_{H^1(\Om)^\prime, H^1(\Om)} \; dt\\
        &=-\int_0^T \into \sum_{0 \le j \ne i \le n} L_{ij} \Big[u_{\eps, j} \nabla u_{\eps, i}- u_{\eps, i} \nabla u_{\eps, j}{- u_{\eps, i} u_{\eps, j} \nabla z_i(\bbu_\eps)}+ {J}_{\eps,ij} \Big] \cdot \nabla \varphi_i \; \d x dt,
		\end{split}
    \end{equation}
	{where $z_{i}$ is given by \eqref{z}.}
	\end{enumerate}
	
\end{definition}
 Our first result is to show global-in-time existence of such solutions. 

\begin{theorem}
    \label{thm:existence-nl}

Let {$\bbu^0_{\eps} \in H^1(\Om)^{n+1}$} be  such that $u_{\eps, i}^0 \ge 0$ for all $i=0, \dots, n$ and $\sum\limits_{i=0}^n u_{\eps,i}^0 = 1$. 
Then, there exists a solution $(\bbu_{\eps}, {J}_{\eps})$ to \eqref{scalar-nonlocal}  in the sense of Definition \ref{def-nonlocal}.
    
\end{theorem}

 For the proof of Theorem \ref{thm:existence-nl} we refer to Section \ref{sect3} below.

\begin{definition}[Solution to the local system]
\label{def-local}

Let  $T>0$ be fixed and let {$\bbu^0 \in H^2(\Om)^{n+1}$} be an initial condition satisfying the bounds of \eqref{in-cond}. We say that $\bbu= (u_{0}, \dots, u_{ n})$ is a solution to the local system \eqref{scalar-local} if

	\begin{enumerate}
	
	\item[$(1_{\text{L}})$] $0 \le u_{ i} \le 1 $ for all $i=0, \dots, n$ and $\displaystyle \sum_{i=0}^n u_{i}= 1$ a.e. in $\Om \times (0,T)$; 
	\item[$(2_{\text{L}})$]  $u_{ i} \in {L^2(0,T; H^2(\Om))}$,  $\partial_t u_{ i} \in L^2(0,T; {H^1(\Omega)'})$ as well as \[
        u_{ i} u_{ j} \nabla\sum_{k=0}^n(c_{ik}-c_{jk})\Delta u_{k} \in L^2(0,T;L^2(\Om))
    \]
    for all $i,j=0, \dots, n$;

	\item[$(3_{\text{L}})$] 
    {$u_{i}(\cdot \,, t)\to u_{i}^0 $
    strongly in $L^2(\Omega)$ as $t\to 0$,}
    for all $i= 0, \dots, n$;
	
	\item[$(4_{\text{L}})$] for every $i= 0, \dots, n $ and every $\varphi_i \in L^2(0,T; H^1(\Om))$ there holds
		\begin{equation*}
		\begin{split}
		& \int_0^T \langle \partial_t u_{i}, \varphi_i \rangle_{H^{1}(\Om)^\prime, H^{1}(\Om)} \; dt \\
		& = -\int_0^T \into \sum_{0 \le j \ne i \le n} L_{ij} \Big[u_{ j} \nabla u_{ i}- u_{ i} \nabla u_{ j}- u_{i} u_{j} \nabla (z_i(\bbu)-z_j(\bbu))+ u_{ i} u_{ j} \nabla\sum_{k=0}^n(c_{ik}-c_{jk})\Delta u_{k}\Big] \cdot \nabla \varphi_i \d x \d t.
		\end{split}
		\end{equation*}
	
	\end{enumerate}
	
\end{definition}

 Our second main result concerns nonlocal-to-local convergence of solutions. 

\begin{theorem}
\label{theorem-local}
 Let $\bbu^0 \in H^2(\Om)$. For $\eps>0$ let $\bbu_\eps$ be a solution to the nonlocal system in the sense of Definition \ref{def-nonlocal} with initial datum $\bbu_\eps^0 = \bbu^0$ and such that \eqref{(a)}--\eqref{e1} are satisfied. Then, up to subsequences, there exists $\bbu \in L^2(0,T; H^1(\Om)^{n+1}) \cap H^1(0,T; (H^1(\Om)')^{n+1})$, {$\bbu(\cdot \,, t)\to \bbu^0 $
    strongly in $L^2(\Omega)$ as $t\to 0$,} such that
    \begin{align}
        u_{\eps,i} \to u_i \quad \text{weakly in } L^2(0,T; H^1(\Om)) \text{ and weakly in } H^1(0,T; H^1(\Om)') \label{(i)} \\
        B_{\eps}(u_{\eps,i}) \to -\Delta u_i \text{ weakly in } {L^2(0,T; {L^2(\Omega)})} \label{ii)}
    \end{align}
and $u_i$ is a solution to the local system in the sense of Definition \ref{def-local}.
\end{theorem}

 The proof of Theorem \ref{theorem-local} is the focus of Section \ref{sect4}.

\section{ Existence of nonlocal solutions}
\label{sect3}

This section is  devoted  to the proof of Theorem \ref{thm:existence-nl} { whose starting point is a, still non-{linear}, regularized time-discrete scheme.}
Existence of iterates is guaranteed by a fixed-point argument involving the chemical potentials {as unknowns}. 
Finally, using the fact {that }the formal gradient flow structure of the system {is preserved in the discrete scheme}, we obtain a priori estimates on the iterates which are sufficient to pass to the limit as $\tau \to 0^+$ and recover the desired weak solution.

Let us emphasize the fact that a key role is played by Lemma \ref{lemma:J}, in which we treat the convergence of the least regular terms of our system. To obtain {the desired} convergence, we will need to use a careful truncation argument along with the strong convergence in $L^2(0,T; L^2(\Om))$ as well as the continuity properties of the operator $B_{\eps}$.

\subsection{Time-discrete approximation}

We start by considering the closed, convex and non-empty subset of $L^{ 2}(\Om)^{n+1}$ given by

	\begin{equation*}
	\mathcal A= \left\{\bbu \in L^{2}(\Om)^{n+1}: \, u_i \ge 0, \ i= 0, \dots, n, \text{ and }  \sum_{i=0}^n u_i = 1\;\text{a.e. in }\Omega \right\}
 \end{equation*}
 endowed with the $L^2$-topology. Note that in particular if $\bbu \in \mathcal{A}$, then for every $i= 0, \dots, n$ it holds $0 \le u_i \le 1$  almost everywhere. 
Further, we introduce below the set of admissible variations, defined as
\begin{equation*}
    \mathcal B= \left\{\boldsymbol{\phi} \in L^{\infty}(\Om)^{n+1}: \, \phi_0=-  \sum_{i=1}^n \phi_i\;\text{a.e.~in }\Omega  \right\}.
 \end{equation*}
 We first prove an  existence result for iterates of a regularized, time-discrete approximation of \eqref{eq:nonlocal}.
 In the statement below, we will denote by $\langle\cdot,\cdot,\rangle_{H^2(\Omega)}$ the scalar product in $H^2(\Omega)$.

\begin{proposition}
    \label{thm:existence-discr}

Let $\tau>0$ be a discrete time step, let $p \in \N$ and let $\bbu_{\eps}^p \in \mathcal{A}$. Then, there exists a solution $(\bbu_{\eps}^{p+1}, \bmu_{\eps}^{p+1}) \in \mathcal{A}  \times H^2(\Om)^{n+1} $ to the following coupled system: for every $ i=0, \dots, n$ and $\varphi_i \in H^2(\Om)$,
    	\begin{align}
	 \into \frac{u_{\eps,i}^{p+1}- u_{\eps, i}^p}{\tau} \varphi_i \; \d x &= - \into \Big(\sum_{0 \le j \ne i \le n} L_{ij} u_{\eps, i}^{p+1} u_{\eps_j}^{p+1} \nabla (\mu_{\eps, i}^{p+1}- \mu_{\eps, j}^{p+1})\Big) \cdot \nabla \varphi_i \; \d x \nonumber  \\
	& \qquad - \tau  \langle \mu_{\eps, i}^{p+1}, \varphi_i \rangle_{H^2(\Om)},  \label{reg}
	\end{align}
and for every $\boldsymbol{\phi} \in \mathcal{B}$,
    \begin{align}\label{eq:mu_to_u}
        \sum_{i=0}^n \into \Big(\ln u_{\eps,i}^{p+1} + \sum_{k=0}^n c_{ik} B_\eps(u_{\eps,k}^{{p+1}})\Big) \phi_i \; \d x= \sum_{i=0}^n \into [\mu_{\eps,i}^{p+1}  {+ \sum_{k=0}^n a_{ik}u_{k}^p}]\phi_i\; \d x.
    \end{align}
Additionally, {$\sum_{i=0}^n\mu_{\eps,i}^{p+1}=0$ a.e.~in $\Omega$ and} there exists $\delta_p>0$ such that
    \[
    u_{\eps,i}^{p+1} \ge \delta_{p+1} \quad \text{for all } i=0, \dots, n, \text{ a.e. in }  \Om.
    \]
\end{proposition}
\begin{remark}
{Note that in \eqref{eq:mu_to_u}, the contribution stemming from the concave part of the energy, $\sum_{k=0}^n a_{ik}u_{k}^p$, appears as given datum on the right-hand side. This implicit-explicit splitting will ensure that the time-discrete energy remains non-increasing.}

{Also note that due to the definition of $\mathcal{B}$, \eqref{eq:mu_to_u} only corresponds to $n$ independent equations, thus only determining $n$ of the $n+1$ functions $u_{\eps,i}$. The remaining one is again given by the constraint $\sum_{i=0}^n u_{\eps,i} = 1$.}
\end{remark}
The idea of the proof is to first linearize \eqref{reg}, { with respect to the unknowns $\mu_{\eps,i}$}, by replacing $\bbu^{p+1}_{ \eps}$ by a given $\tilde \bbu \in \mathcal{A}$. This gives rise to an operator 
$$
\mathcal S_1\colon \mathcal A \subset L^{ 2}(\Om)^{n+1} \to H^2(\Om)^{n+1},\quad \mathcal S_1(\tilde \bbu) = \bmu,
$$
{which we show to be well-defined.} 
Next, using variational methods, we establish that for every given $\bmu$ there exists a solution $\bbu$ to \eqref{eq:mu_to_u} defining an operator
$$
\mathcal{S}_2 \colon {L^\infty}(\Om)^{n+1} \to \mathcal{A},\quad \mathcal{S}_2(\bmu) = \bbu.
$$
A fixed point of the composed operator $\mathcal S_2 \circ{\imath_{H^2\to L^\infty}\circ} \mathcal S_1$, {with $\imath_{H^2\to L^\infty} \colon u\in H^2(\Omega)\mapsto u\in L^\infty(\Omega)$ being the canonical compact Sobolev immersion operator,} is  then  a weak solution to \eqref{reg} and its existence will be a consequence of Schauder's fixed point theorem. 

We will structure the proof in several lemmas below, starting with the linearized version of \eqref{reg}.

\begin{lemma}
\label{lemma:mu}

Let $\widetilde \bbu \in \mathcal A$. Then, there exists a unique solution $\bmu \in H^2(\Om)^{n+1}$ to 
		\begin{equation}
		\label{ex-1}
		\begin{split}
		\into \frac{\widetilde u_i- u_{ \eps,i}^p}{\tau} \varphi_i \; \d x&= -\into \Big(\sum_{0 \le j \ne i \le n} L_{ij} \widetilde u_i \widetilde u_j \nabla (\mu_i- \mu_j) \Big) \cdot \nabla \varphi_i \; \d x- \tau \langle \mu_i, \varphi_i \rangle_{H^2(\Om)}
		\end{split}
		\end{equation}
for every $i= 0, \dots, n$ and every $\varphi_i \in H^2(\Om)$.  Moreover, {$\sum_{i=0}^n\mu_i=0$ a.e.~in $\Omega$ and} there exists $C_1= C_1(n, \tau, \Om)>0$  such that
	\begin{equation}
	\label{ex-1.1}
	\|\bmu\|_{(H^2(\Om))^{n+1}} \le C_1.
	\end{equation}

\end{lemma}

\begin{proof}

Thanks to \eqref{mobility}, system \eqref{ex-1} can be written in vectorial form as follows
	\begin{equation}
	\label{ex-2}
	-\frac{1}{\tau} \into (\widetilde \bbu- \bbu^p_{\eps}) \cdot \varphi \; \d x= \into \nabla \varphi \cdot  {M}(\widetilde \bbu) \nabla \bmu \; \d x+ \tau \langle \varphi, \bmu \rangle_{H^2(\Om)^{n+1}},
	\end{equation}
for  every  $\varphi \in H^2(\Om)^{n+1}$. Moreover,
	\begin{equation}
	\label{ex-3}
	0 \le  {M}(\widetilde \bbu) \le (n+1) \, L \,  {{I}}_{n+1} \quad \text{a.e. in } \Om
	\end{equation}
in the sense of symmetric matrices, where $L= \max\limits_{0 \le j \ne i \le n} L_{ij}$ and ${{I}}_{n+1}$ is the identity matrix of $\R^{(n+1) \times (n+1)}$. Then, existence and uniqueness to \eqref{ex-2} is a consequence of the Lax-Milgram theorem, \cite{Brezis}. 
{The fact that $\sum_{i=0}^n\mu_{\eps,i}^{p+1}=0$ a.e.~in $\Omega$ follows easily by choosing $\varphi_i=\varphi$ independent of $i$ in \eqref{ex-2}, summing in $i$ and exploiting the definition of $\A$ as well as the property $L_{ij}=L_{ji}$. Finally, choosing}
$\varphi= \bmu$ in \eqref{ex-2} and using the semidefiniteness of $\bM$ as well as H$\ddot{\text{o}}$lder inequality gives
	\[
	\tau \|\bmu\|_{H^2(\Om)^{n+1}}^2 \le \frac{1}{\tau} \sum_{i=0}^n \|\widetilde u_i- u_{\eps,i}^p\|_{L^2(\Om)} \|\mu_i\|_{L^2(\Om)} \le \frac{ C}{\tau} \Big(\sum_{i=0}^n \|\widetilde u_i- u_{\eps,i}^p\|_{L^2(\Om)}\Big) \|\bmu\|_{H^2(\Om)^{n+1}}.
	\]
Since $\widetilde \bbu, \bbu^p \in \mathcal A$ it follows that
	\[
	\|\bmu\|_{(H^2(\Om))^{n+1}} \le \frac{ C}{\tau^2} 2 (n+1)|\Om|^{1/2},
	\]
and therefore \eqref{ex-1.1} is satisfied.
\end{proof}

This result shows that the operator $\mathcal S_1\colon \mathcal A \subset L^\infty(\Om)^{n+1} \to H^2(\Om)^{n+1}$ which associates to any $\widetilde \bbu \in \mathcal A$ the unique solution $\bmu$ to \eqref{ex-1} is well-defined. The following result guarantees its continuity and as it follows almost verbatim as in \cite[Lemma 3.5]{EMP}, we omit its proof here.  Note that the set $\mathcal A$ here is endowed with the $L^2$-topology whereas in \cite[Lemma 3.5]{EMP} it is considered as a subset of $L^{\infty}(\Om)^{n+1}$. Nevertheless, {arguing as in the proof of } \cite[Lemma 3.5]{EMP} yields Lemma \ref{lem:s1_continuous} below.

\begin{lemma}\label{lem:s1_continuous}
The map $\mathcal S_1\colon \mathcal A \to H^2(\Om)^{n+1}$ is continuous.
\end{lemma}

As next step we now aim to recover $\bbu \in \mathcal A$ from a given $ \bmu \in {L^\infty}(\Om)^{n+1}$ by means of \eqref{eq:mu_to_u}.
We will tackle this problem by using a variational approach.  To be precise,  we will identify  $\bbu\in \A$  as the unique solution
to the minimization problem
	\begin{equation}
	\label{ex-6}
	\min_{\bbw \in \mathcal A} F_{\bmu}(\bbw),
	\end{equation}
where $F_{\bmu}\colon \mathcal A \to \R $ is defined by
	\[
	F_{\bmu} (\bbw)= \sum_{i=0}^n \into w_i \ln w_i +\frac{1}{2}\sum_{k=0}^n c_{ik} B_{\eps}(w_k) w_i - [\mu_i {+ \sum_{k=0}^n a_{ik}u_{k}^p}]w_i \; \d x.
	\]
We first show that this problem indeed admits a unique minimum.
\begin{lemma}
\label{lemma:ex-un-u}
    Let $\bmu \in {L^\infty}(\Om)^{n+1}$. Then, there exists a unique solution $\bbu$ to \eqref{ex-6}.
\end{lemma}

{Ok for me}
\begin{proof}
We use the direct methods of the Calculus of Variations. 
First of all, given {$\bmu \in L^\infty(\Omega)^{n+1}$} we show that $F_{\bmu}$ is bounded from below. Indeed, we use the definition of $\A$, the fact that $x \ln x \ge 1/e$ in $[0,1]$, a Young's inequality in the nonlocal terms and $\bmu \in  L^\infty(\Om)^{n+1}$ to infer 
    \[
    F_{\bmu}(\bbw) \ge  \frac{n}{e} |\Om|- 2n \max_{0 \le i,k \le n} |c_{ik}|  \|K_\eps\|_{L^1(\Om^2)}- \left(\sum_{i=0}^n \|\mu_i\|_{L^\infty(\Om)} 
    {+}\max_{0 \le i,k \le n}|a_{ik}|\right) |\Om|> -\infty.
    \]
Observe further that the function 
$\bbw^*=(1/(n+1), \ldots, 1/(n+1))$ 
is an element of $\A$ and satisfies 
$F_{\bmu}(\bbw^*) < +\infty$. Thus, it holds that 
    \[
    -\infty< \inf_{\mathcal{A}} F_{\bmu} \le F_{\bmu} (\bbw^*)< +\infty.
    \]
{Therefore $F_{\bmu}$ is bounded from below and not identically $+\infty$. We show now that $F_{\bmu}$ is strictly convex. Clearly the term $\sum_{i=0}^n\int_\Omega w_i\log w_i \d x$ is strictly convex. On the other hand the nonlocal term can be rewritten as}
\begin{align}
\nonumber
\sum_{i,k=0}^n \into c_{ik} B_{\eps} (u_k) u_i \; \d x&= \sum_{i,k=0}^n \into c_{ik} ((K_{\eps} \ast 1) u_k- K_{\eps} \ast u_k) u_i \; \d x \\
\nonumber
&= \frac{1}{2}\sum_{i,k=0}^n \into \into c_{ik} K_{\eps}(x,y) (u_k(x)- u_k(y)) (u_i(x)- u_i(y)) \; \d x\d y \\
\label{Beps.quadform}
&= \frac{1}{2}\into \into K_{\eps} (x,y) (\bbu(x)- \bbu(y))^t {\cc} (\bbu(x)- \bbu(y)) \; \d x\d y,
\end{align}
and in turn the positive semidefiniteness of the matrix 
${\cc}$ implies the convexity of the above term. 
{Finally, the last term is linear in $w$, hence convex. Therefore $F_{\bmu}$ is strictly convex. The continuity of $F_{\bmu}$ with respect to the strong topology of $L^2(\Omega)^{n+1}$ follows easily from the fact that $\bmu, \bbu^p\in L^2(\Omega)^{n+1}$ and from the above-written representation of the nonlocal term. It follows that $F_{\bmu}$ is lower semicontinuous for the weak topology in $L^2(\Omega)^{n+1}$. From standard arguments in the Direct Method of the Calculus of Variations we infer that $F_{\bmu}$ admits a minimizer $\bbu\in\mathcal{A}$. The minimizer is unique thanks to the strict convexity of $F_{\bmu}$. This finishes the proof of the Lemma.}
\end{proof}

In order to identify \eqref{eq:mu_to_u} as the Euler-Lagrange equation to \eqref{ex-6}, and to  obtain  that minimizers of $F_{\bmu}$ are weak solutions to \eqref{eq:mu_to_u}, we need to be able to construct perturbations in such a way that we do not leave the set $\mathcal A$. This will be possible due to pointwise bounds of the minimizers which only depend on $\bmu$ through its $L^\infty$-norm as the following result shows.

\begin{lemma}\label{lem:bounds_minimizer}

For  every  $\bmu \in {L^\infty}(\Om)^{n+1}$ there exists $\delta_{\bmu}>0$ such that for any minimizer $\bbu $ to $F_{\bmu}$ in $\A$ there  holds
    \[
    u_i \ge \delta_{\bmu} \quad \text{for all } i= 0, \dots, n  \text{ and  almost everywhere  in } \Om.
    \]
 Additionally,  for all $N>0$, there exists $\delta>0$ which only depends on $\displaystyle N, n, \max_{0 \le i,j \le n} |c_{ij}|, \max_{0 \le i,j \le n} |a_{ij}|, \|K_{\eps}\|_{L^1(\Om^2)}$, such that for all $\bmu \in {L^\infty}(\Om)^{n+1}$ with $\|\bmu\|_{{L^\infty}(\Om)^{n+1}} \le N$ and for any minimizer $\bbu$ to $F_{\bmu}$  there holds 
    \[
    u_i \ge \delta \quad \text{for all } i=0, \dots, n \text{ and  almost everywhere  in } \Om.
    \]
    
\end{lemma}

\begin{proof}

 Let $\bmu  \in L^{\infty}(\Om)^{n+1}$ and let $\bbu $ be a minimizer of $F_{\bmu}$ on $\mathcal A$. 
We aim to show that there exists $\delta >0$ 
such that $\delta \leq u_i$ almost everywhere in $\Omega$ for all $i=0, \dots, n$. The precise dependence of $\delta$ on the data will be specified later in the proof. \\
To fix the ideas, we will show the result for $u_0$, being understood that the same reasoning applies to the other components. 
We reason by contradiction and assume that the Lebesgue measure of the set
$
\mathcal M_\delta := \{ x \in \Omega \; : \; u_0(x) < \delta \}
$ is positive. 
Now, let us define 
\begin{align}  %
u_0^\delta := \max (u_0, \delta) \quad \text{as well as} \quad
\begin{split}\label{eq:u1u2delta}
u_i^\delta := u_i - (u_0^\delta - u_0)\frac{u_i}{1-u_0},\quad i=1,\ldots, n,
\end{split}
\end{align}
and set $\bbu^\delta := (u_0^\delta, \ldots, u_n^\delta)$. In (\ref{eq:u1u2delta}), since $1-u_0 = \sum\limits_{j=1}^n u_j\geq u_i \geq 0$, the function $\frac{u_i}{1-u_0}$ is well-defined almost everywhere using the convention that 
$\frac{u_i}{1-u_0} = 0$ as soon as $u_i = 0$. By definition,  we have  $1 \geq u_0^\delta \geq 0$ and $u_0^\delta + \sum\limits_{i=1}^nu_i^\delta = 1$. 
Furthermore,  $u_i^\delta(x) =0$ for all $x\in \Omega$ such that $u_i(x) = 0$. For all $x\in \Omega$ such that $u_i(x)>0$,  it follows that $1-u_0(x) \geq u_i(x)>0$ and 
$$
	u_i^\delta(x) = u_i(x)\left(1 - \frac{u_0^\delta(x) - u_0(x)}{1-u_0(x)}\right) \ge 0, \quad \text{ since } \quad \frac{u_0^\delta(x) - u_0(x)}{1-u_0(x)} \le \frac{1 - u_0(x)}{1-u_0(x)} = 1.
	$$
As a consequence,  $\bbu^\delta \in \mathcal A$ and $u_0^\delta = 1 - \sum\limits_{i=1}^n u_i^\delta$. 
We now prove that for $\delta$ sufficiently small, $F_ {\bmu}(\bbu^\delta) < F_{\bmu}(\bbu)$. 
Using the fact 
that $u_i^\delta = u_i$ for $i=0,\ldots,n$ on $\mathcal M_\delta^c = \{x \in \Omega\; : \; u_0(x) \ge \delta \}$ yields 
	\begin{align}
	F_{\bmu}(\bbu^\delta) - F_{\bmu}(\bbu) 
	= &\sum_{i=0}^n \int_{\mathcal M_\delta}\left( 
    u_i^\delta \ln u_i^\delta - u_i \ln u_i 
    \right)\d x
    -\sum_{i=0}^n \int_{\mathcal M_\delta}
    \left(\mu_i - \sum_{k=0}^n a_{ik}u_{k}^p\right) (u_i^\delta-u_i)\d x
    \nonumber\\ 
    &+\int_\Omega\sum_{i,k=0}^n c_{ik} \big(B_{\eps}(u_k^{\delta})u_i^{\delta}- B_{\eps}(u_k)u_i \big)\d x
    \label{eq:terms_linear}
	\end{align}
For the nonlocal terms, we write
\begin{align*}
&\sum_{i,k=0}^n \int_{\Omega} c_{ik} \big(B_{\eps}(u_k^{\delta}) u_i^{\delta}- B_{\eps}(u_k) u_i \big) \; \d x \\
&=\sum_{i,k=0}^n \int_{\Omega} c_{ik} 
B_{\eps}(u_k^{\delta}-u_k) u_i^{\delta} \; \d x +
\sum_{i,k=0}^n \int_{\Omega} c_{ik} 
B_{\eps}(u_k) (u_i^{\delta}-u_i) \; \d x
\end{align*}
Given \eqref{Beps.esteps} and the bounds $0\leq u_i^\delta\leq u_i\leq 1$ we infer
\begin{align*}
&\sum_{i,k=0}^n \int_{\Omega} c_{ik} \big(B_{\eps}(u_k^{\delta}) u_i^{\delta}- B_{\eps}(u_k) u_i \big) \; \d x\leq C\int_{\Omega}|\bbu^\delta - \bbu|\d x
=C\int_{\mathcal{M}_\delta}|\bbu^\delta - \bbu|\d x
\leq C\int_{\mathcal{M}_\delta}|u_0^\delta - u_0|\d x.
\end{align*}
The last terms in \eqref{eq:terms_linear} are easily estimated by means of \eqref{eq:u1u2delta}, {using $[\mu_i - \sum_{k=0}^n a_{ik}u_{k}^p] \in L^\infty(\Omega)$}, while for the logarithmic terms we follow the proof of \cite[Lemma 3.7]{EMP} so that 
\eqref{eq:terms_linear} reduces to
    \[
    F_{\bmu}(\bbu^{\delta})- F_{\bmu}(\bbu) \le (\ln \delta+ C) \int_{\mathcal{M}_{\delta}} (u_0^{\delta}- u_0) \; \d x, 
    \]
where $C>0$ is a constant depending on the data but not on $\delta$. Moreover, note that $\displaystyle \int_{\mathcal{M}_{\delta}} (u_0^{\delta}- u_0) \; \d x>0$, because the function $u_0^{\delta}- u_0>0$ on $\mathcal{M}_{\delta}$, which has positive measure. Therefore, if $\delta$ is chosen small enough to guarantee that $\ln \delta+ C<0$, it follows that
    \[
    F_{\bmu}(\bbu^{\delta})- F_{\bmu}(\bbu)< 0,
    \]
which contradicts the fact that $\bbu$ is a minimizer to $F_{\bmu}$. 

Since the same argument applies to all the remaining components, we can conclude that for every minimizer $\bbu \in \mathcal{A}$ to $F_{\bmu}$ there exists $\delta_{\bmu}>0$ such that
    \[
    u_i \ge \delta_{\bmu} \quad \text{for all } i= 0, \dots, n. 
    \]
Arguing as in the proof of \cite[Lemma 3.7]{EMP}, end of Step 1, we see that $\delta$ can be chosen to be dependent only on the data, as soon as $\bmu$ is assumed to satisfy $\|\bmu\|_{{L^\infty}(\Om)^{n+1}} \le N$.
This concludes the proof of the lemma.
\end{proof}

\begin{lemma}

Let $\bbu \in \mathcal{A} $ be the unique minimizer of $F_{\bmu}$. Then, $\bbu$ is a weak solution to \eqref{eq:mu_to_u}.
In particular, 
\begin{align}
\label{sol2}
\mu_i = 
\ln u_i+ \sum_{k=0}^n c_{ik} B_\eps(u_k)- \sum_{k=0}^n a_{ik}u_{k}^p + \lambda
\quad \text{a.e. in } \Om\text{ for all $i=0, \dots, n$},
\end{align}
where $\lambda = -\frac{1}{n+1}\sum_{i=0}^n(\ln u_i+ \sum_{k=0}^n c_{ik} B_\eps(u_k)- \sum_{k=0}^n a_{ik}u_{k}^p)$ is a Lagrange multiplier ensuring that the Gauge condition $\sum_{i=0}^n\mu_i = 0$ is fulfilled.
\end{lemma}

\begin{proof}

Since $\bbu$ minimizes $F_{\bmu}$ { in $\mathcal{A}$},  property \eqref{eq:mu_to_u} follows arguing exactly as in the proof of \cite[Lemma 3.8]{EMP}. 
{Now fix $i\in\{1,\ldots,n\}$, $\psi\in C^\infty_c(\Omega)$ arbitrary. Choosing $\bphi\in\mathcal{B}$ such that $\bphi_j = \delta_{ij}\psi$ for $j\in\{1,\ldots,n\}$ and $\bphi_0=-\psi$ and employing the fundamental Lemma of the Calculus of Variations yield \eqref{sol2}.}
\end{proof}

{
\begin{lemma}
If $\bmu, \bbu^p\in H^1(\Omega)^{n+1}$ then every minimizer $\bbu$ to $F_{\bmu}$ in $\mathcal{A}$ satisfies $\bbu\in H^1(\Omega)^{n+1}$.
\end{lemma}
\begin{proof}
Let $D_h$ be a standard difference quotient operator, i.e.~$D_h f(x) \equiv |h|^{-1}(f(x+ h)-f(x))$ for every function $f : \Omega\to\R$ and every vector $h\in\R^d\backslash\{0\}$. 
Choosing $\phi_i = D_{-h}D_h u_{\eps,i}$ in \eqref{eq:mu_to_u} and performing a discrete integration by parts yields
\begin{align*}
\int_\Omega\sum_{i=0}^n D_h\ln u_i\, D_h u_i \d x
+ \int_\Omega\sum_{i,k=0}^n c_{ik}D_h u_i\cdot B_\eps(D_h u_k)\d x 
= \int_\Omega\sum_{i=0}^n D_h u_i\cdot\left(D_h\mu_i + \sum_{k=0}^n a_{ik}D_h u_k^p\right)\d x
\end{align*}
The second integral on the left-hand side of the above inequality is nonnegative thanks to the positive definiteness of $\cc$ and \eqref{B.qf}, while the first one can be bound from below via the elementary inequality $(x-y)\log(x/y)\geq (x-y)^2$, which holds for every $x,y\in (0,1]$:
$$
\int_\Omega\sum_{i=0}^n D_h\ln u_i\, D_h u_i \d x\geq 
\int_\Omega \sum_{i=0}^n |D_h u_i|^2\d x.
$$
Therefore via Cauchy-Schwartz inequality we infer
$$
\int_\Omega \sum_{i=0}^n |D_h u_i|^2\d x\leq C \int_\Omega\sum_{i=0}^n 
(|D_h u_i^p|^2 + |D_h\mu_i|^2)\d x.
$$
Standard results on Sobolev spaces yield the statement. This finishes the proof of the Lemma.
\end{proof}
}

Let now $\mathcal{S}_2 \colon {L^\infty}(\Om)^{n+1} \to \mathcal{A}$ be the map that associates each $\bmu \in {L^\infty}(\Om)^{n+1}$ to the unique minimizer of $F_{\bmu}$ in $\mathcal{A}$, that is, the unique weak solution $\bbu \in \mathcal{A}$ to \eqref{eq:mu_to_u}. We have the following regularity result.

\begin{lemma}\label{lem:s2_continuous}
The operator $\mathcal{S}_2 \colon {L^\infty}(\Om)^{n+1} \to \mathcal{A}$ is continuous. 
    
\end{lemma}

\begin{proof}

Let $(\bmu^m)_{m \in \N} \subset {L^\infty}(\Om)^{n+1}$ be such that $\bmu^m \to \bmu$ strongly in ${L^\infty}(\Om)^{n+1}$. Set $\bbu^m= \mathcal{S}_2(\bmu^m)$ as well as $\bbu= \mathcal{S}_2(\bmu)$. We want to show that $\mathcal{S}_2(\bmu^m) \to \mathcal{S}_2(\bmu)$  strongly in $L^2(\Om)^{n+1}$  as $m \to \infty$. 
{
Writing \eqref{eq:mu_to_u} for $\bbu^m$ and $\bbu$, taking the difference between the two equations and choosing $\bphi=\bbu^m-\bbu$ in the resulting identity leads to
} 
\begin{equation*}
\sum_{i=0}^n\int_{\Omega}(\ln u_i^m-\ln u_i)(u^m_i-u_i)\; \d x
{+\int_\Omega\sum_{i,j=0}^n c_{ij}B_\eps(u_i^m-u_i)(u_j^m-u_j)\d x}    
\leq 2\into \sum_{i=0}^n|\mu^m_i-\mu_i|\; \d x,
\end{equation*}
{where we bound the right-hand side by employing the fact that $0\leq u_i, u_i^m\leq 1$ a.e.~in $\Omega$. The second integral in the above identity is nonnegative thanks to \eqref{Beps.quadform}, hence we focus on the first one. The elementary inequality $(x-y)\log(x/y)\geq ({x}-{y})^2$, which holds for every $x,y\in (0,1]$, together with the strong $L^1(\Omega)$ convergence $\bmu^m\to\bmu$ imply $\bbu^m\to\bbu$ strongly in $L^2(\Omega)^{n+1}$. This finishes the proof of the Lemma.
}
\end{proof}
We now have all the necessary ingredients at hand to prove the existence of iterates.
\begin{proof}[Proof of Proposition \ref{thm:existence-discr}]
By the results of Lemmas \ref{lemma:mu}, \ref{lem:s1_continuous}, \ref{lemma:ex-un-u} and \ref{lem:s2_continuous}, together with the compact embedding of $H^2(\Omega)$ into ${L^\infty}(\Omega)$ for $d \le 3$, it follows that the operator $\mathcal S := \mathcal S_2 \circ{\imath_{H^2\to L^\infty}\circ} \mathcal S_1 : \A \to \A$
is both compact and continuous and thus, by Schauder's fixed point theorem, it admits a fixed point $\bbu_{\eps}^{p+1}$ which is the solution to \eqref{reg}.
\end{proof}

\subsection{A priori estimates}

In view of Proposition \ref{thm:existence-discr}, we deduce that for any initial datum $\bbu^0_{\eps}\in \A$ and any $\tau\in (0,1)$ there exists a sequence $(\bbu_{\eps}^p,\mathbb{\mu}_{\eps}^p)_p\subset \A\times H^2(\Om)^{n+1}$, defined recursively as solutions to \eqref{reg} and \eqref{eq:mu_to_u} for every $p\in \N$. 

We then define several piecewise-constant-in-time functions as follows: for all $p\in \mathbb{N} \setminus \{0\}$, for all $i= 0, \dots, n$ and all $t \in ((p-1)\tau, p\tau]$, we set\\
\begin{equation}
\label{setting}
\begin{split}
\bbu_{\eps}^{(\tau)}(t)  &= \bbu_{\eps}^p,\quad  u_{\eps,i}^{(\tau)}(t)  = u_{\eps,i}^p,\\
\mu_{\eps,i}^{(\tau)} =  \mu_{\eps,i}^p&= \ln u_{\eps,i}^p {- \sum_{k=0}^n a_{ik}u_{\eps,k}^{p-1}}+ \sum_{j=0}^n c_{ij} \big((K_{\eps} \ast 1) u_{\eps, j}^p- K_{\eps} \ast u_{\eps,j}^p\big) + \lambda_\eps^p\\
&= \ln u_{\eps,i}^p {- \sum_{k=0}^n a_{ik}u_{\eps,k}^{p-1}} + \sum_{j=0}^n c_{ij} B_{\eps} (u_{\eps,j}^p)+ \lambda_\eps^p,
\end{split}
\end{equation}
for some Lagrange multiplier $\lambda_\eps^p$ which ensures that the condition $\sum_{i=0}^n\mu_{\eps,i}^p = 0$ is satisfied.

We also define $\bbu^{(\tau)}(0) =\bbu^0$. Let $P^{(\tau)} = \lfloor T/\tau \rfloor \in \mathbb{N} \setminus \{0\}$ be the largest integer such that $P^{(\tau)}\tau\leq T$. 
We additionally introduce the time-shifted solution $\sigma_\tau \bbu_{\eps}^{(\tau)}$ as
$$
\sigma_\tau \bbu_{\eps}^{(\tau)}(t) = \bbu_{\eps}^{p-1} \mbox{ for all } t \in ((p-1)\tau, p\tau],  \, p\in \mathbb{N} \setminus \{0\},
$$
whose components are given by  $(\sigma_\tau u_{\eps,0}^{(\tau)}, \dots, \sigma_\tau u_{\eps,n}^{(\tau)})$. For all $\bbu_{\eps} = (u_{\eps,0}, u_{\eps,1}, \dots, u_{\eps,n}) \in \mathcal A$, for all $\tau >0$ and {$t>\tau$} we define the entropy functional
	\begin{align*}
	E_{\text{NL}}^{(\tau)}(t):= &\sum_{i=0}^n \into u_{\eps,i}^{(\tau)}(t) \ln u_{\eps,i}^{(\tau)}(t) \; \d x\\ 
    &+ \frac{1}{4} \sum_{i,j=0}^n \into \into c_{ij} K_{\eps}(x,y) (u_{\eps,i}^{(\tau)}(t,x)- u_{\eps,i}^{(\tau)}(t,y)) (u_{\eps,j}^{(\tau)}(t,x)- u_{\eps,j}^{(\tau)}(t,y)) \; \d x\d y \\
    & -\frac{1}{2}\sum_{i,k=0}^n\int_\Omega a_{ik} \sigma_\tau u_{\eps,i}^{(\tau)}\sigma_\tau u_{\eps,k}^{(\tau)}\;\d x =:E_{\text{NL}}^{(\tau),\text{conv}}(t) + E_{\text{NL}}^{(\tau),\text{conc}}(t), 
	\end{align*}
where $E_{\text{NL}}^{(\tau),\text{conv}}(t)$ collects the first two convex terms. We also define
	\begin{align*}
	E_{\text{NL}}^{p+1}&= \sum_{i=0}^n \into u_{\eps,i}^{p+1} \ln u_{\eps,i}^{p+1} \; \d x+ \frac{1}{4} \sum_{i,j=0}^n \into \into c_{ij} K_{\eps}(x,y) (u_{\eps,i}^{p+1}(x)- u_{\eps,i}^{p+1}(y)) (u_{\eps,j}^{p+1}(x)- u_{\eps,j}^{p+1}(y)) \; \d x\d y\\
    &{-\frac{1}{2}\sum_{i,k=0}^n\int_\Omega a_{ik}  u_{\eps,i}^{p}u_{\eps,k}^{p}}\;\d x {=:E_{\text{NL}}^{\text{conv},p+1} + E_{\text{NL}}^{\text{conc},p}}. 
	\end{align*}
{such that for all $p\in \mathbb{N}$ we have  $E_{\text{NL}}^{(\tau)}(t_{p+1}) = E_{\text{NL}}^{p+1}$.}
The main contribution of this subsection is to establish some a priori bounds on the interpolants $(u^{(\tau)}_{\eps,i},{\mu}^{(\tau)}_{\eps,i})$, with constants independent of $\tau$ and $\eps$. 
\begin{remark}
    \label{rem:energy}
Note that the energy functional $E_{\text{NL}}$ in \eqref{E} is bounded from below, since the matrix $\cc$ is positive semidefinite and the function $x \mapsto x\ln x- x+ 1$ is bounded from below in $[0,1]$.  
\end{remark}

We begin this subsection by stating the monotonicity of $ E_{\text{NL}}$.

\begin{lemma}
\label{lemma:monot}

Let $({\bbu}_{\eps}^{p})_{p \in \N}$ be a sequence of solutions to \eqref{reg}. Then, the sequence  $ (E_{\text{NL}}^{(\tau)}(t_p))_{p\in\N}$ 
is decreasing. Moreover, there exists a constant $C> 0$ such that for $i=0,\ldots, n$
	\[
	\into \into K_{\eps}(x, y) (u_{\eps, i}^{p+1}(x)- u_{\eps, i}^{p+1}(y))^2 dx \d y \le C.
	\]

\end{lemma}

\begin{proof}
We reason as in \cite[Lemma 4.3]{EMP}.
We test each equation of \eqref{reg} with $\varphi_i= \mu_{\eps, i}^{p+1}$ and then sum for $i= 0, \dots, n$ to have
	\begin{align}
	\sum_{i=0}^n \into \frac{u_{\eps,i}^{p+1}- u_{\eps, i}^p}{\tau} \mu_{\eps, i}^{p+1} \; \d x=& 
    {- \sum_{i=0}^n \into \Big( \sum_{0 \le j \ne i \le n} L_{ij} u_{\eps, i}^{p+1} u_{\eps, j}^{p+1} \nabla(\mu_{\eps, i}^{p+1}- \mu_{\eps, j}^{p+1}) \Big) \cdot \nabla \mu_{\eps, i}^{p+1} \; \d x}\nonumber \\
	& - \tau \sum_{i=0}^n   \|\mu_{\eps, i}^{p+1}\|^2_{H^2(\Om)}. 	\label{estim1}
	\end{align}
For the right-hand side, we exploit the fact the mobility $\bM$ is positive semidefinite, see  \cite[Remark 4.1]{EMP}, to have %
	\begin{align}
	& - \sum_{i=0}^n \into \Big( \sum_{0 \le j \ne i \le n} L_{ij} u_{\eps, i}^{p+1} u_{\eps, j}^{p+1} \nabla(\mu_{\eps, i}^{p+1}- \mu_{\eps, j}^{p+1}) \Big) \cdot \nabla \mu_{\eps, i}^{p+1} \; \d x 
    - \tau \sum_{i=0}^n   \|\mu_{\eps, i}^{p+1}\|^2_{H^2(\Om)}\nonumber \\
	& = -\int_\Omega (\nabla {\bmu}_{\eps}^{p+1})^T {\bM}({\bbu}_{\eps}^{p+1}) \nabla {\bmu}_{\eps}^{p+1}\;\d x 
    - \tau \sum_{i=0}^n   \|\mu_{\eps, i}^{p+1}\|^2_{H^2(\Om)}\le 0. 	\label{estim2}
	\end{align}
On the left-hand side we use \eqref{mu-i}, {using the convexity and concavity of the respective parts of the energy}, we obtain
	\begin{equation}
	\label{estim3}
	\begin{split}
	&\sum_{i=0}^n \into \frac{u_{\eps, i}^{p+1}-u_{\eps, i}^{p}}{\tau}  \mu_{\eps, i}^{p+1} \; \d x\\
    &= \frac{1}{\tau} \sum_{i=0}^n \into (u_{\eps, i}^{p+1}-u_{\eps, i}^{p}) D_{u_{\eps, i}^{p+1}} E_{\text{NL}}^{\text{conv},p+1} -\frac{1}{\tau} \sum_{i=0}^n \into (u_{\eps, i}^{p}-u_{\eps, i}^{p+1}) D_{u_{\eps, i}^{p}} E_{\text{NL}}^{\text{conc},p} \; \d x  \\
    &\ge \frac{1}{\tau} \big(E_{\text{NL}}^{p+1}- E_{\text{NL}}^p\big).
	\end{split}
	\end{equation}
From \eqref{estim1}-\eqref{estim3} we infer 
	\[
	\frac{1}{\tau} \big(E_{\text{NL}}^{p+1}- E_{\text{NL}}^p\big)\le 0,
	\]
that is, the sequence $(E_{\text{NL}}(\bbu_\eps^p))_{p \in \N}$ is non-increasing. In particular, there exists $C_{NL}>0$ such that $E_{\text{NL}}(\bbu_\eps^{p+1}) \le E_{\text{NL}}(\bbu_\eps^0) \le C_{NL}$ which, taking \eqref{E} {and \eqref{lb.C}} into account, entails
	\[
	\frac{{\eta_C}}{4}{\sum_{i=0}^n} \into \into K_{\eps}(x, y) (u_{\eps, i}^{p+1}(x)- u_{\eps, i}^{p+1}(y))^2 dx \d y \le C_{NL}.
	\]
This completes the proof.
\end{proof}

We now exploit this property to establish some a priori estimates which play a central role in order to pass to the limit as $\tau \to 0^+$. The first result is a direct consequence of the fact that the discrete energy is decreasing.
For the sake of the next proposition we define the quantity
\begin{equation}
\label{eq:gradient-u0}
\mathcal{K}_\eps^0 :=
\sum_{i=0}^n \into \into K_\eps(x,y) |\nabla u_{\eps,i}^{0}(x)- \nabla u_{\eps,i}^{0}(y)|^2 \; \d x \d y.
\end{equation}
Notice that, or every $\eps>0$, $\mathcal{K}_\eps^0<\infty$ since $u_{\eps,i}^{0}\in H^1(\Omega)$.
\begin{proposition}
\label{thm:estim}
Let $({\bfu_{\eps}^{p}})_{p \in \N}$ be a sequence of solutions to \eqref{reg}. Then, there exists a constant $C>0$, depending on $\mathcal{K}_\eps^0$, $\|\bbu_\eps^0\|_{H^1(\Om)}$, $\Omega$, $n$, $A$, $L$ and $\cc$ such that
 \begin{align}
	&\sum_{i=0}^n \int_0^T \into \frac{|\nabla u_{\eps,i}^{(\tau)}|^2}{u_{\eps,i}^{(\tau)}} \; \d x dt \le C, \label{a} \\
	&\sum_{i=0}^n \int_0^T \into \into  K_{\eps}(x, y) |\nabla u_{\eps,i}^{(\tau)}(x)- \nabla u_{\eps,i}^{(\tau)}(y)|^2 \; \d x \d y dt \le C, \label{b} \\
 	& {\sum_{i=0}^n \sum_{0 \le j \ne i \le n} \int_0^T \into |u_{\eps, i}^{(\tau)} u_{\eps, j}^{(\tau)} \nabla (q_{\eps,i}(\bbu_{\eps}^{(\tau)})-  q_{\eps,j}(\bbu_{\eps}^{(\tau)}))|^2 \; \d x dt \le C,} \label{e} \\
	& \tau \sum_{i=0}^n \int_0^T \| \mu_{\eps, i}^{(\tau)}\|_{H^2(\Om)}^2 \; dt \le C. \label{d}
	\end{align}
{Furthermore, if the quantity $\mathcal{K}_\eps^0$ in \eqref{eq:gradient-u0} and $\|\bbu_\eps^0\|_{H^1(\Om)}$ are uniformly bounded in $\eps\in (0,1)$, then the constant $C$ in \eqref{a}--\eqref{d} can be chosen to be independent of $\eps$.}
\end{proposition}

\begin{proof}

We use an argument already exploited in \cite{EMP}. 
{We start from \eqref{estim1}.}
The left-hand side {of \eqref{estim1}} can be {estimated as in \eqref{estim3}.}
Concerning the right-hand side,  we first set $\displaystyle \ell:= \min_{\overset{0 \le i, j \le n}{i \ne j}} L_{ij}$. We have
	\begin{align}
	\label{e2}
    & -\sum_{i=0}^n \into \sum_{0 \le j \ne i \le n} L_{ij} u_{\eps, i}^{p+1} u_{\eps, j}^{p+1} \nabla (\mu_{\eps, i}^{p+1}- \mu_{\eps, j}^{p+1}) \cdot \nabla \mu_{\eps, i}^{p+1} \; \d x - \tau \sum_{i=0}^n \|\mu_{\eps,i}^{p+1}\|_{H^2(\Om)}^2 \\
	& = - \sum_{i=0}^n \into \sum_{0 \le j \ne i \le n} (L_{ij}- \ell) u_{\eps, i}^{p+1} u_{\eps, j}^{p+1} \nabla (\mu_{\eps, i}^{p+1}- \mu_{\eps, j}^{p+1}) \cdot \nabla \mu_{\eps, i}^{p+1} \; \d x \nonumber \\
	& \qquad - \ell \sum_{i=0}^n \into \sum_{0 \le j \ne i \le n} u_{\eps, i}^{p+1} u_{\eps, j}^{p+1} \nabla (\mu_{\eps, i}^{p+1}- \mu_{\eps, j}^{p+1}) \cdot \nabla \mu_{\eps, i}^{p+1} \; \d x -\tau \sum_{i=0}^n \|\mu_{\eps, i}^{p+1}\|_{H^2(\Om)}^2. \nonumber
	\end{align}
The first term reads as
	\begin{align}
	\label{e3}
	&- \sum_{i=0}^n \into \sum_{0 \le j \ne i \le n} (L_{ij}- \ell) u_{\eps, i}^{p+1} u_{\eps, j}^{p+1} \nabla (\mu_{\eps, i}^{p+1}- \mu_{\eps, j}^{p+1}) \cdot \nabla \mu_{\eps, i}^{p+1} \; \d x \\
 &= - \into (\nabla {\bmu_{\eps}^{p+1}})^T \widetilde {\bM}(\bbu_{\eps}^{p+1}) \nabla {\bmu}_{\eps}^{p+1} \le 0, \nonumber
	\end{align}
where $\widetilde {\bM}(\bbu_{\eps}^{p+1}) $ is given as in \eqref{mobility}, but with the coefficients $L_{ij}$ replaced by $L_{ij}- \ell$. 
{By using \eqref{setting}, \eqref{q}, \eqref{z}, we rewrite} the second term reads as follows
\begin{align}
\label{e4}
-& \ell \sum_{i=0}^n \into \sum_{0 \le j \ne i \le n} u_{\eps, i}^{p+1} u_{\eps, j}^{p+1} \nabla (\mu_{\eps, i}^{p+1}- \mu_{\eps, j}^{p+1}) \cdot \nabla \mu_{\eps, i}^{p+1} \; \d x \\
=& - \ell \sum_{i=0}^n \into \sum_{0 \le j \ne i \le n} u_{\eps,i}^{p+1} u_{\eps, j}^{p+1} \nabla (\ln u_{\eps, i}^{p+1}- \ln u_{\eps, j}^{p+1}) \cdot \nabla \ln u_{\eps,i}^{p+1} \; \d x \nonumber \\
& - \ell \sum_{i=0}^n \into \sum_{0 \le j \ne i \le n} u_{\eps,i}^{p+1} u_{\eps, j}^{p+1} \nabla (\ln u_{\eps, i}^{p+1}- \ln u_{\eps, j}^{p+1}) \cdot \nabla q_{\eps,i}(\bbu^{p+1}) \; \d x \nonumber \\
&- \ell \sum_{i=0}^n \into \sum_{0 \le j \ne i \le n} u_{\eps,i}^{p+1} u_{\eps, j}^{p+1}  \nabla (q_{\eps,i}(\bbu^{p+1})- q_{\eps,j}(\bbu^{p+1}))\cdot \nabla \ln u_{\eps, i}^{p+1}  \; \d x \nonumber \\
&- \ell \sum_{i=0}^n \into \sum_{0 \le j \ne i \le n} 
u_{\eps,i}^{p+1} u_{\eps, j}^{p+1} \nabla (q_{\eps,i}(\bbu^{p+1})- q_{\eps,j}(\bbu^{p+1})) \cdot \nabla q_{\eps,i}(\bbu^{p+1}) \; \d x \nonumber\\
&+\ell \sum_{i=0}^n \into \sum_{0 \le j \ne i \le n} u_{\eps,i}^{p+1} u_{\eps, j}^{p+1} \nabla (\ln u_{\eps, i}^{p+1}- \ln u_{\eps, j}^{p+1}) \cdot \nabla z_i(\bbu^{p}) \; \d x\nonumber \\
&+\ell \sum_{i=0}^n \into \sum_{0 \le j \ne i \le n} u_{\eps,i}^{p+1} u_{\eps, j}^{p+1} \nabla (z_i(\bbu^p) - z_j(\bbu^p)) \cdot \nabla (\ln(u_{\eps, i}^{p+1}) - z_i(\bbu^p)) \; \d x\nonumber \\
&+\ell \sum_{i=0}^n \into \sum_{0 \le j \ne i \le n} u_{\eps,i}^{p+1} u_{\eps, j}^{p+1} \nabla (z_i(\bbu^p) - z_j(\bbu^p) \cdot \nabla  q_{\eps,i}(\bbu^{p+1})) \; \d x\nonumber \\
&+\ell \sum_{i=0}^n \into \sum_{0 \le j \ne i \le n} u_{\eps,i}^{p+1} u_{\eps, j}^{p+1}  \nabla (q_{\eps,i}(\bbu^{p+1})- q_{\eps,j}(\bbu^{p+1}))\cdot \nabla z_i(\bbu^p)  \; \d x 
\nonumber \\ 
=:&\mathrm{DE}_1 + \mathrm{DE}_2 + \mathrm{DE}_3+ \mathrm{DE}_4+ \mathrm{DE}_5+ \mathrm{DE}_6+ \mathrm{DE}_7+ \mathrm{DE}_8
\end{align}
We are going to estimate the several terms in the above equation separately. {We start by those terms which come from the convex parts of the energy.} First of all, 
{in the arguments that follow we make use, without any further mention of it, of the elementary identity $\sum_{0 \le j \ne i \le n} u_j = 1-u_i$ (an immediate consequence of \eqref{constraint}). }
{Next, let us compute for $i=0,\ldots,n$
\begin{align}
\sum_{0 \le j \ne i \le n} u_{\eps,i}^{p+1} u_{\eps, j}^{p+1} \nabla (\ln u_{\eps, i}^{p+1}- \ln u_{\eps, j}^{p+1}) =& 
\sum_{0 \le j \ne i \le n} (u_{\eps,j}^{p+1} \nabla u_{\eps, i}^{p+1}-  u_{\eps, i}^{p+1} \nabla u_{\eps, {j}}^{p+1})\nonumber\\
=& (1-u_{\eps,i}^{p+1}) \nabla u_{\eps, {i}}^{p+1}+  u_{\eps, i}^{p+1} \nabla u_{\eps, i}^{p+1}
= \nabla u_{\eps, i}^{p+1}.
\label{new.nalnu}
\end{align}
From \eqref{new.nalnu} we immediately deduce
}
\begin{align*}
\mathrm{DE}_1
= & - \ell \sum_{i=0}^n \into \sum_{0 \le j \ne i \le n} u_{\eps,i}^{p+1} u_{\eps, j}^{p+1} \nabla (\ln u_{\eps, i}^{p+1}- \ln u_{\eps, j}^{p+1}) \cdot \nabla \ln u_{\eps,i}^{p+1} \; \d x 
=  -\ell \sum_{i=0}^n \into \frac{|\nabla u_{\eps,i}^{p+1}|^2}{u_{\eps,i}^{p+1}} \; \d x.
\end{align*}
For the second term use again $\eqref{new.nalnu}$ 
{and employ \eqref{B.qf}, \eqref{lb.C} to deduce}
\begin{align*}
\mathrm{DE}_2 =& -\ell \sum_{i=0}^n \into \sum_{0 \le j \ne i \le n}  u_{\eps,i}^{p+1} u_{\eps,j}^{p+1} \nabla(\ln u_{\eps,i}^{p+1}- \ln u_{\eps,j}^{p+1}) \cdot \nabla q_{\eps,i}(\bbu_\eps^{p+1}) \; \d x \\
=& -\ell \sum_{i=0}^n \into \nabla u_{\eps,i}^{p+1} \cdot \nabla \Big(\sum_{k=0}^n c_{ik} B_{\eps}(u_{\eps,k}^{p+1})\Big)= -{\ell}\sum_{i,k=0}^n \into c_{ik} \nabla u_{\eps,i}^{p+1} \cdot B_\eps(\nabla u_{\eps,k}^{p+1}) \; \d x \\
=& {- \frac{\ell}{2} \sum_{i,k=0}^n \into \into  c_{ik} K_\eps(x,y) \nabla (u_{\eps,i}^{p+1}(x)- u_{\eps,i}^{p+1}(y)) \cdot \nabla(u_{\eps,k}^{p+1}(x)- u_{\eps,k}^{p+1}(y)) \; \d x \d y}\\
\le &-{\frac{\ell\eta_C}{2}} \sum_{i=0}^n \into \into K_\eps(x,y) |\nabla u_{\eps,i}^{p+1}(x)- \nabla u_{\eps,i}^{p+1}(y)|^2 \; \d x \d y.
\end{align*}
For the third term, {we apply a Young's inequality weighted with $(u_{\eps, i}^{p+1})^{1/2}$} to have
    \begin{align*}
    \mathrm{DE}_3 =& - \ell \sum_{i=0}^n \into \sum_{0 \le j \ne i \le n} u_{\eps,i}^{p+1} u_{\eps, j}^{p+1} \nabla \ln u_{\eps, i}^{p+1} \cdot \nabla (q_{\eps,i}(\bbu_{\eps}^{p+1})- q_{\eps,j}(\bbu_{\eps}^{p+1})) \; \d x \\
    =& - \ell \sum_{i=0}^n \into \sum_{0 \le j \ne i \le n}  u_{\eps, j}^{p+1} \nabla  u_{\eps, i}^{p+1} \cdot \nabla (q_{\eps,i}(\bbu_{\eps}^{p+1})- q_{\eps,j}(\bbu_{\eps}^{p+1})) \; \d x \\
    \le &{\frac{\ell}{3}} \sum_{i=0}^n \sum_{0 \le j \ne i \le n} \into u_{\eps,i}^{p+1} u_{\eps,j}^{p+1} |\nabla (q_{\eps,i}(\bbu_{\eps}^{p+1})-  q_{\eps,j}(\bbu_{\eps}^{p+1}))|^2 \; \d x
    + {\frac{3\ell}{4}} \sum_{i=0}^n \sum_{0 \le j \ne i \le n} \into u_{\eps, j}^{p+1} \frac{|\nabla u_{\eps,i}^{p+1}|^2}{u_{\eps, i}^{p+1}} \; \d x \\
    \le &{\frac{\ell}{3}} \sum_{i=0}^n \sum_{0 \le j \ne i \le n} \into u_{\eps,i}^{p+1} u_{\eps,j}^{p+1} |\nabla (q_{\eps,i}(\bbu_{\eps}^{p+1})-  q_{\eps,j}(\bbu_{\eps}^{p+1}))|^2 \; \d x
    + {\frac{3\ell}{4}} \sum_{i=0}^n  \into  \frac{|\nabla u_{\eps,i}^{p+1}|^2}{u_{\eps, i}^{p+1}} \; \d x.
    \end{align*}
    
Exploiting symmetries in the fourth term gives
    \begin{align*}
     \mathrm{DE}_4 =& - \ell \sum_{i=0}^n \into \sum_{0 \le j \ne i \le n} u_{\eps,i}^{p+1} u_{\eps, j}^{p+1} \nabla (q_{\eps,i}(\bbu_{\eps}^{p+1})- q_{\eps,j}(\bbu_{\eps}^{p+1})) \cdot \nabla q_{\eps,i}(\bbu_{\eps}^{p+1}) \; \d x \\
     =& - {\frac{\ell}{2}} \sum_{i=0}^n \sum_{0 \le j \ne i \le n} \into u_{\eps, i}^{p+1} u_{\eps, j}^{p+1} |\nabla \big(q_{\eps,i}(\bbu_{\eps}^{p+1})- q_{\eps,j}(\bbu_{\eps}^{p+1})\big)|^2 \; \d x.   
    \end{align*}
{We now continue to estimate terms that contain parts from the concave part of the energy, i.e. contributions of the previous time step $p$. We will repeatedly make use of inequality \eqref{eq:h1_control_nonlocal}. }
For the fifth and sixth term, we estimate {via \eqref{new.nalnu}}
\begin{align*}
    \mathrm{DE}_5  &= \ell \sum_{i=0}^n \into \sum_{0 \le j \ne i \le n} u_{\eps,i}^{p+1} u_{\eps, j}^{p+1} \nabla (\ln u_{\eps, i}^{p+1}- \ln u_{\eps, j}^{p+1}) \cdot \nabla z_i(\bbu^{p}) \; \d x \\
    &\le C \left(\sum_{i=0}^n \|\nabla u_{\eps,i}^{p}\|_{L^2(\Om)}^2 + \sum_{i=0}^n \|\nabla u_{\eps,i}^{p+1}\|^2_{L^2(\Om)}\right),\\ %
\mathrm{DE}_6 &=\ell \sum_{i=0}^n \into \sum_{0 \le j \ne i \le n} u_{\eps,i}^{p+1} u_{\eps, j}^{p+1} \nabla (z_i(\bbu^p) - z_j(\bbu^p)) \cdot \nabla (\ln(u_{\eps, i}^{p+1}) - z_i(\bbu^p)) \; \d x\\
&\le C \left(\sum_{i=0}^n \|\nabla u_{\eps,i}^{p}\|_{L^2(\Om)}^2 + \sum_{i=0}^n \|\nabla u_{\eps,i}^{p+1}\|^2_{L^2(\Om)}\right). %
\end{align*}
Next, we apply symmetry arguments to obtain
\begin{align*}
\mathrm{DE}_7&=\ell \sum_{i=0}^n \into \sum_{0 \le j \ne i \le n} u_{\eps,i}^{p+1} u_{\eps, j}^{p+1} \nabla (z_i(\bbu^p) - z_j(\bbu^p)) \cdot \nabla q_{\eps,i}(\bbu^{p+1})) \; \d x \\
&=\frac{\ell}{2} \sum_{i=0}^n \into \sum_{0 \le j \ne i \le n} u_{\eps,i}^{p+1} u_{\eps, j}^{p+1} \nabla (z_i(\bbu^p) - z_j(\bbu^p)) \cdot (\nabla q_{\eps,i}(\bbu^{p+1})-\nabla q_{\eps,j}(\bbu^{p+1}))
\; \d x\\
&\le C \sum_{i=0}^n \|\nabla u_{i,\eps}^p\|_{L^2(\Om)}^2 + \frac{\ell}{{12}}\sum_{i=0}^n \sum_{0 \le j \ne i \le n} \into u_{\eps,i}^{p+1} u_{\eps,j}^{p+1} |\nabla (q_{\eps,i}(\bbu_{\eps}^{p+1})-  q_{\eps,j}(\bbu_{\eps}^{p+1}))|^2 \; \d x.
\end{align*}
For the last term we estimate
\begin{align*}
\mathrm{DE}_8&=\ell \sum_{i=0}^n \into \sum_{0 \le j \ne i \le n} u_{\eps,i}^{p+1} u_{\eps, j}^{p+1}  \nabla (q_{\eps,i}(\bbu^{p+1})- q_{\eps,j}(\bbu^{p+1}))\cdot \nabla z_i(\bbu^p)  \; \d x \\
&\le \frac{\ell}{{24}}\sum_{i=0}^n \sum_{0 \le j \ne i \le n} \into u_{\eps,i}^{p+1} u_{\eps,j}^{p+1} |\nabla (q_{\eps,i}(\bbu_{\eps}^{p+1})-  q_{\eps,j}(\bbu_{\eps}^{p+1}))|^2 \; \d x + C\sum_{i=0}^n\|\nabla u_i^p\|_{L^2(\Om)}^2.
\end{align*}
Thus, recalling again \eqref{estim3}, we obtain
    \begin{align*}
       & \frac{E_{{NL}}(\bbu_\eps^{p+1})- E_{{NL}}(\bbu_\eps^p)}{\tau} 
       + {\frac{\ell\eta_C}{2}}
       \sum_{i=0}^n \into \into  K_\eps(x,y) |\nabla u_{\eps,i}^{p+1}(x,s)- \nabla u_{\eps,i}^{p+1}(y,s)|^2 \; \d x \d y \\
       & \quad +\frac{\ell}{{4}} \sum_{i=0}^n \into \frac{|\nabla u_{\eps,i}^{p+1}|^2}{u_{\eps,i}^{p+1}} \; \d x
       + \frac{\ell}{{24}} \sum_{i=0}^n \sum_{0 \le j \ne i \le n} \into u_{\eps, i}^{p+1} u_{\eps, j}^{p+1} |\nabla \big(q_{\eps,i}(\bbu_{\eps}^{p+1})- q_{\eps,j}(\bbu_{\eps}^{p+1})\big)|^2 \; \d x \\
        &\quad +\tau \sum_{i=0}^n \|\mu_{\eps, i}^{p+1}\|_{H^2(\Om)}^2 \le  C \left(\sum_{i=0}^n \|\nabla u_{\eps,i}^{p}\|_{L^2(\Om)}^2 + \sum_{i=0}^n \|\nabla u_{\eps,i}^{p+1}\|^2_{L^2(\Om)}\right).
    \end{align*}
We continue by using \eqref{eq:h1_control_nonlocal} {with $\delta = \frac{\ell\eta_C}{4C}$ and $\delta = \frac{\ell\eta_C}{8C}$ for the term involving $u_{\eps,i}^{p+1}$ and $u_{\eps,i}^{p}$, respectively,} to further estimate the right-hand side as follows
\begin{align*}
C &\left(\sum_{i=0}^n \|\nabla u_{\eps,i}^{p}\|_{L^2(\Om)}^2 + \sum_{i=0}^n \|\nabla u_{\eps,i}^{p+1}\|^2_{L^2(\Om)}\right)\\
\le &{\frac{\ell\eta_C}{4}}\sum_{i=0}^n \into \into K_\eps(x,y) |\nabla u_{\eps,i}^{p+1}(x)- \nabla u_{\eps,i}^{p+1}(y)|^2 \; \d x \d y 
+ C \|u_{\eps,i}^{p+1}\|_{L^2(\Om)}^2\\
&+{\frac{\ell\eta_C}{8}}\sum_{i=0}^n \into \into K_\eps(x,y) |\nabla u_{\eps,i}^{p}(x)- \nabla u_{\eps,i}^{p}(y)|^2 \; \d x \d y 
+ C \|u_{\eps,i}^{p}\|_{L^2(\Om)}^2,
\end{align*}
hence {by putting the two above inequalities together and using the fact that $0\leq u_{\eps,i}^p, u_{\eps,i}^{p+1}\leq 1$ a.e.~in $\Omega$} we obtain
\begin{align*}
\frac{1}{\tau}& \big(E_{{NL}}(\bbu_\eps^{p+1})- E_{{NL}}(\bbu_\eps^p)\big) \\
& +\frac{\ell}{{4}} \sum_{i=0}^n \into \frac{|\nabla u_{\eps,i}^{p+1}|^2}{u_{\eps,i}^{p+1}} \; \d x +{\frac{\ell\eta_C}{4}} \sum_{i=0}^n \into \into K_\eps(x,y) |\nabla u_{\eps,i}^{p+1}(x)- \nabla u_{\eps,i}^{p+1}(y)|^2 \; \d x \d y \\
& + {\frac{\ell}{24}} 
\sum_{i=0}^n \sum_{0 \le j \ne i \le n} \into u_{\eps, i}^{p+1} u_{\eps, j}^{p+1} |\nabla \big(q_{\eps,i}(\bbu_{\eps}^{p+1})- q_{\eps,j}(\bbu_{\eps}^{p+1})\big)|^2 \; \d x +\tau \sum_{i=0}^n \|\mu_{\eps, i}^{p+1}\|_{H^2(\Om)}^2\\
\le & C + {\frac{\ell\eta_C}{8}}\sum_{i=0}^n \into \into K_\eps(x,y) |\nabla u_{\eps,i}^{p}(x)- \nabla u_{\eps,i}^{p}(y)|^2 \; \d x \d y.
\end{align*}
We then multiply this {inequality} by $\tau$, sum for $0 \le p \le P^{(\tau)}-1$ and use {the fact that the energy functional $E_{NL}$ is bounded from below (see Remark \ref{rem:energy})} to obtain 
    \[
    \begin{split}
        & \frac{\ell}{{4}} \sum_{i=0}^n \int_0^T \into \frac{|\nabla u_{\eps,i}^{(\tau)}|^2}{u_{\eps,i}^{(\tau)}} \; \d x dt+ {\frac{\ell\eta_C}{8}} \sum_{i=0}^n \int_0^T \into \into K_\eps(x,y) |\nabla u_{\eps,i}^{(\tau)}(x)- \nabla u_{\eps,i}^{(\tau)}(y)|^2 \; \d x \d y dt \\
        & \quad + \frac{\ell}{{24}} \sum_{i=0}^n \sum_{0 \le j \ne i \le n} \int_0^T \into u_{\eps, i}^{(\tau)} u_{\eps, j}^{(\tau)} |\nabla \big(q_{\eps,i}(\bbu_{\eps}^{(\tau)})- q_{\eps,j}(\bbu_{\eps}^{(\tau)})\big)|^2 \; \d x dt +\tau \sum_{i=0}^n \int_0^T \|\mu_{\eps, i}^{(\tau)}\|_{H^2(\Om)}^2 \; dt \\
        & \le C(T+1) + E_{{NL}}(\bbu_{\eps}^0) + {\frac{\ell\eta_C}{8}}\sum_{i=0}^n \into \into K_\eps(x,y) |\nabla u_{\eps,i}^{0}(x)- \nabla u_{\eps,i}^{0}(y)|^2 \; \d x \d y\\
        &\le C(T+1) {+\frac{\ell\eta_C}{8}\mathcal{K}_\eps^0 + C\|\bbu_{\eps}^0\|_{H^1(\Omega)}^2}, %
    \end{split}
    \]
where {in the last passage} we used \eqref{eq:gradient-u0} {and inequality $E_{{NL}}(\bbu_{\eps}^0)\leq C(1+\|\bbu_{\eps}^0\|_{H^1(\Omega)}^2)$, which follows from \eqref{energy-NL} and \cite[Thm.~1]{BBM}.}
Since all the quantities on the left-hand side are nonnegative, the bounds \eqref{a}, \eqref{b} and \eqref{d} are satisfied, while the bound \eqref{e} follows by noticing that, since $u_{\eps, i}^{(\tau)}, u_{\eps,j}^{(\tau)} \in (0,1)$, then
    \[
   \int_0^T \into |u_{\eps, i}^{(\tau)} u_{\eps, j}^{(\tau)} \nabla \big(q_{\eps,i}(\bbu_{\eps}^{(\tau)})- q_{\eps,j}(\bbu_{\eps}^{(\tau)})\big)|^2 \; \d x dt \le \int_0^T \into u_{\eps, i}^{(\tau)} u_{\eps, j}^{(\tau)} |\nabla \big(q_{\eps,i}(\bbu_{\eps}^{(\tau)})- q_{\eps,j}(\bbu_{\eps}^{(\tau)})\big)|^2 \; \d x dt,
    \]
for all $i=0, \dots, n$ and all $j \ne i= 0, \dots, n$.
This completes the proof.
\end{proof}

\begin{lemma}
\label{lemma:dt-tau}
    There exists $C>0$, independent of $\tau>0 $ and $\eps>0$, such that
        \[
        \int_0^T \Big\|\frac{u_{\eps,i}^{(\tau)}- \sigma_\tau u_{\eps,i}^{(\tau)}}{\tau}\Big\|_{H^2(\Om)^\prime}^2 \; dt \le C \quad \text{for all } i=0, \dots, n.
        \]
\end{lemma}

\begin{proof}

Follows using the same arguments as in \cite[Lemma 4.7]{EMP}.
\end{proof}

\subsection{The limit as $\tau \to 0^+$}

In this subsection we will recover a weak solution to \eqref{scalar-nonlocal} in the sense of Definition \ref{def-nonlocal} as the weak limit of some extracted subsequence $(\bbu_\eps^{(\tau)})_{\tau>0}$ as $\tau \to 0^+$.

After interpolation, \eqref{reg} reads as
    	\begin{align}
	\int_0^T &\into \frac{u_{\eps,i}^{(\tau)}- \sigma_\tau u_{\eps, i}^{(\tau)}}{\tau} \varphi_i \; \d x dt \nonumber\\
	=& - \int_0^T \into \sum_{0 \le j \ne i \le n} L_{ij} \Big[u_{\eps, j}^{(\tau)} \nabla u_{\eps, i}^{(\tau)}- u_{\eps, i}^{(\tau)} \nabla u_{\eps, j}^{(\tau)}+ u_{\eps, i}^{(\tau)} u_{\eps, j}^{(\tau)} \nabla \Big(\sum_{k=0}^n (c_{ik}-c_{jk}) B_{\eps}(u_{\eps,k}^{(\tau)})\Big)\Big] \cdot \nabla \varphi_i \; \d x dt \nonumber \\
	& +\int_0^T\int_\Om 
    \sum_{0 \le j \ne i \le n} L_{ij}
    u_{\eps, i}^{(\tau)} u_{\eps, j}^{(\tau)} \nabla \left(\sum_{k=0}^n (a_{ik}-a_{jk}) (\sigma_\tau u_{\eps,k}^{(\tau)})\right)\cdot \nabla \varphi_i\;\d x\d t - \tau \int_0^T \langle \mu_{\eps, i}^{(\tau)}, \varphi_i \rangle_{H^2(\Om)} \; dt, \label{discr-interp}
	\end{align}
for every $\varphi_i \in L^2(0,T; H^2(\Om))$.
The estimates collected in the previous section imply the existence of a function {$\bbu_\eps \in L^2(0,T; H^1(\Om))^{n+1}$} such that, up to the extraction of a subsequence, 
    \begin{align}
    \label{tau-uep}
    & u_{\eps,i}^{(\tau)} \rightharpoonup u_{\eps, i} \quad \text{weakly in } L^2(0,T; H^1(\Om)) { \quad \text{as well as weak-}* \text{ in } L^\infty(0,T; L^\infty(\Om)),} \\
    \label{tau-dt-uep}
    & {\frac{u_{\eps,i}^{(\tau)}- \sigma_\tau u_{\eps, i}^{(\tau)}}{\tau} \rightharpoonup \partial_t u_{\eps,i} \quad { \text{weakly in } L^2(0,T; H^2(\Om)^\prime)}}.
    \end{align}
Taking the compact embedding $H^1 (\Om) \hookrightarrow L^2(\Om)$ into account, \cite[Theorem 1]{DJ} implies that 
    \begin{equation}
        \label{lim_ui}
        u_{\eps,i}^{(\tau)} \to u_{\eps,i} \quad \text{strongly in } L^2(0,T; L^2(\Om)),
    \end{equation}
    while the uniform bound of $(u_{\eps,i}^{(\tau)})_{\tau>0}$ in $L^{\infty}(0,T; L^{\infty}(\Om))$ implies that $u_{\eps,i}^{(\tau)} \to u_{\eps,i}$ strongly in $L^r(0,T; L^r(\Om))$, for any $r \in [1, \infty)$ and  for any $i=0, \dots, n$. {The same holds true for $\sigma_\tau u_{\eps,i}^{(\tau)}$, as can be shown by restricting to $t>\tau$ and using the fact that the corresponding indicator function is uniformly bounded.}
From \eqref{Beps.esteps} and the strong convergence $u_{\eps,i}^{(\tau)} \to u_{\eps,i}$ in $L^r(0,T; L^r(\Om))$ for every $r \in [1, \infty)$ we infer that 
$B_\eps(u_{\eps,i}^{(\tau)}) \to B_\eps(u_{\eps,i})$ strongly in $L^r(0,T; L^r(\Om))$ for every $r \in [1, \infty)$. Given \eqref{q}, this also implies that
\begin{align}
\label{q.cnv.tau}
q_{\eps,i}(\bbu_{\eps}^{(\tau)}) \to q_{\eps,i}(\bbu_{\eps})\quad\mbox{strongly in }L^r(0,T; L^r(\Om)),~~\forall r \in [1, \infty).
\end{align}
From \eqref{d} we have that $\big(\sqrt \tau \mu_{\eps,i}^{(\tau)}\big)_{\tau>0}$ is bounded in $L^2(0,T; H^2(\Om))$, therefore
\begin{equation}
\label{TAU}
\tau \mu_{\eps,i}^{(\tau)} \to 0 \quad \text{strongly in } L^2(0,T; H^2(\Om)).
\end{equation}
These convergences enable us to pass to the limit  in \eqref{discr-interp} as $\tau \to 0$. Indeed, for all $i= 0, \dots, n$ and as $\tau \to 0$, up to the extraction of subsequences we have %
    \begin{align}
        & \int_0^T \into \frac{u_{\eps,i}^{(\tau)}- \sigma_{\tau} u_{\eps, i}^{(\tau)}}{\tau} \varphi_i \; \d x dt \to \int_0^T \langle \partial_t u_{\eps,i}, \varphi_i \rangle_{H^2(\Om)', H^2(\Om)} \; dt, \nonumber \\
        & \int_0^T \into \sum_{0 \le j \ne i \le n} L_{ij} \left[u_{\eps, j}^{(\tau)} \nabla u_{\eps, i}^{(\tau)}- u_{\eps, i}^{(\tau)} \nabla u_{\eps, j}^{(\tau)}
        -u_{\eps, i}^{(\tau)} u_{\eps, j}^{(\tau)} \nabla \left(\sum_{k=0}^n (a_{ik}-a_{jk}) (\sigma_\tau u_{\eps,k}^{(\tau)})\right)
        \right] \cdot \nabla \varphi_i \; \d x \d t 
        \nonumber \\
        & \quad\to \int_0^T \into \sum_{0 \le j \ne i \le n} L_{ij} \left[u_{\eps, j} \nabla u_{\eps, i}- u_{\eps, i} \nabla u_{\eps, j} 
        -u_{\eps, i} u_{\eps, j} \nabla \left(\sum_{k=0}^n (a_{ik}-a_{jk}) u_{\eps,k}\right)
        \right] \cdot \nabla \varphi_i \; \d x \d t ,
        \nonumber \\
        & \tau \int_0^T \langle \mu_{\eps, i}^{(\tau)}, \varphi_i \rangle_{H^2(\Om)} \; dt \to 0. \label{easy-limits}
    \end{align}
\medskip
Passing to the limit in the term 
    \begin{align*}
   & \int_0^T \into \sum_{0 \le j \ne i \le n} L_{ij} u_{\eps, i}^{(\tau)} u_{\eps, j}^{(\tau)} \nabla \Big(\sum_{k=0}^n c_{ik} B_{\eps}(u_{\eps,k}^{(\tau)})- \sum_{l=0}^n c_{jl} B_\eps(u_{\eps,l}^{(\tau)})\Big) \cdot \nabla \varphi_i \; \d x dt \\
   & \qquad \quad = \int_0^T \into \sum_{0 \le j \ne i \le n} L_{ij} u_{\eps, i}^{(\tau)} u_{\eps, j}^{(\tau)} \nabla \Big(q_{\eps,i}(\bbu_{\eps}^{(\tau)})- q_{\eps,j}(\bbu_{\eps}^{(\tau)}   \big)\Big) \cdot \nabla \varphi_i \; \d x dt
    \end{align*}
is more delicate and will be treated in details in the following lemma.

\begin{lemma}
\label{lemma:J}
    For all $i \ne j= 0, \dots, n$ there exists ${J}_{\eps, ij} \in L^2(0,T; L^2(\Om)^d)$ such that $ J_{\eps,ij}= u_{\eps,i} u_{\eps,j} \nabla (q_{\eps,i}(\bbu_\eps)- q_{\eps,j}(\bbu_\eps))$ holds in the weak sense, that is,
    \[
    \int_0^T \into {J}_{\eps,ij} \cdot {\eta} \; \d xdt= - \int_0^T \into (q_{\eps,i}(\bbu_\eps)- q_{\eps,j}(\bbu_\eps)) \divergenz(u_{\eps,i} u_{\eps,j} {\eta}) \; \d x \d t
    \]
for every ${\eta} \in L^2(0,T; H^1(\Om)^d) \cap L^\infty ((0,T) \times \Om; \R^d)$.
{Furthermore it holds (up to the extraction of a subsequence)}
    \[
    u_{\eps, i}^{(\tau)} u_{\eps,j}^{(\tau)} \nabla (q_{\eps,i}(\bbu_\eps^{(\tau)})- q_{\eps,j}(\bbu_\eps^{(\tau)})) \to {J}_{\eps,ij} \quad \text{weakly in } L^2(0,T; L^2(\Om)^d)
    \]
    as $\tau \to 0$.
\end{lemma}

\begin{proof}

Fix $i,j \in \{0, \dots, n \}$ such that $j \ne i $ and set, for simplicity,    
    \begin{equation}
        \label{lim1}
    {J}_{\eps,ij}^{(\tau)}:= u_{\eps,i}^{(\tau)} u_{\eps,j}^{(\tau)} \nabla (q_{\eps,i}(\bbu_\eps^{(\tau)})- q_{\eps,j}(\bbu_\eps^{(\tau)})).
    \end{equation}
Estimate \eqref{e} implies the existence of $ J_{\eps, ij} \in L^2(0,T; L^2(\Om)^d) $ such that
    \begin{equation}
        \label{lim3}
        {J}_{\eps,ij}^{(\tau)} \rightharpoonup  J_{\eps, ij} \quad \text{weakly in } L^2(0,T; L^2(\Om)^d).
    \end{equation}
From \eqref{q.cnv.tau} we infer that $q_{\eps,i}(\bbu_\eps^{(\tau)})- q_{\eps,j}(\bbu_\eps^{(\tau)})\to q_{\eps,i}(\bbu_\eps)- q_{\eps,j}(\bbu_\eps)$ strongly in $L^2(0,T; L^2(\Omega))$, while \eqref{tau-uep} and \eqref{lim_ui} imply that $u_{\eps,i}^{(\tau)}u_{\eps,j}^{(\tau)}\rightharpoonup u_{\eps,i}u_{\eps,j}$ weakly in $L^2(0,T; H^1(\Omega))$. We infer
\begin{align*}
&\int_0^T\int_\Omega {J}_{\eps,ij}\cdot\eta\;\d x\d t = 
\lim_{\tau\to 0}
\int_0^T\int_\Omega {J}_{\eps,ij}^{(\tau)}\cdot\eta\;\d x\d t \\
&= 
-\lim_{\tau\to 0}\int_0^T\int_\Omega (q_{\eps,i}(\bbu_\eps^{(\tau)})- q_{\eps,j}(\bbu_\eps^{(\tau)}))\divergenz (u_{\eps,i}^{(\tau)}u_{\eps,j}^{(\tau)}\eta)\;\d x\d t
\\
&= 
-\int_0^T\int_\Omega (q_{\eps,i}(\bbu_\eps)- q_{\eps,j}(\bbu_\eps))\divergenz (u_{\eps,i}u_{\eps,j}\eta)\;\d x\d t
\end{align*}
for every $\eta\in C^\infty(\Omega)$. This finishes the proof of the Lemma.
\end{proof}

\begin{proof}[Proof of Theorem \ref{thm:existence-nl}]

In order to obtain the desired result, it remains to pass to the limit as $\tau \to 0$ in \eqref{discr-interp} using the convergences in \eqref{easy-limits} and the result of Lemma \ref{lemma:J}.

It follows that there exists $\bbu_{\eps} \in L^2(0,T; H^1(\Om))^{n+1}$ with $\partial_t \bbu_{\eps} \in L^2(0,T; H^2(\Om)')^{n+1}$ such that %
\begin{align*}
\int_0^T &\langle \partial_t u_{\eps,i}, \varphi_i \rangle_{H^2(\Om)', H^2(\Om)} \; dt \nonumber\\
=& - \int_0^T \into \sum_{0 \le j \ne i \le n} L_{ij} \Big[u_{\eps, j} \nabla u_{\eps, i}- u_{\eps, i} \nabla u_{\eps, j}+ u_{\eps, i} u_{\eps, j} \nabla \Big(\sum_{k=0}^n (c_{ik}-c_{jk}) B_{\eps}(u_{\eps,k})\Big)\Big] \cdot \nabla \varphi_i \; \d x dt, \nonumber \\
&+\int_0^T\int_\Om \sum_{0 \le j \ne i \le n} L_{ij} u_{\eps, i} u_{\eps, j} \nabla \left(\sum_{k=0}^n (a_{ik}-a_{jk}) u_{\eps,k}\right)\cdot \nabla \varphi_i\;\d x\d t
\end{align*}
for all $i=0, \dots, n$ and for all $\varphi_i \in L^2(0,T; H^2(\Om))$. The previous weak formulation implies that $\partial_t u_{\eps, i} \in L^2(0,T; H^1(\Om)')$ and therefore by density we can extend the whole formulation to any $\varphi_i \in L^2(0,T; H^1(\Om))$, {thus replacing the duality product $\langle \partial_t u_{\eps,i}, \varphi_i \rangle_{H^2(\Om)', H^2(\Om)}$ with $\langle \partial_t u_{\eps,i}, \varphi_i \rangle_{H^1(\Om)', H^1(\Om)} $.}
{
Given that $\bbu_{\eps} \in L^2(0,T; H^1(\Om))^{n+1}$ and $\partial_t \bbu_{\eps} \in L^2(0,T; H^1(\Om)')^{n+1}$ it follows \cite[Prop.~23.23]{Zei12} that $\bbu_{\eps}\in C^0([0,T]; L^2(\Omega))$. By arguing like in \cite{juengel} one shows that $\bbu_{\eps}(t)\to\bbu^0$ strongly in $H^1(\Omega)'$ as $t\to 0^+$, but since $\bbu_{\eps}\in C^0([0,T]; L^2(\Omega))$ the convergence $\bbu_{\eps}(t)\to\bbu^0$ is strong in $L^2(\Omega)$.
}
Therefore, it turns out that $(\bbu_{\eps}, {J}_{\eps})$ is a weak solution to \eqref{scalar-nonlocal} in the sense of Definition \ref{def-nonlocal}. The proof is thus complete. 
\end{proof}

\section{From nonlocal to local}
\label{sect4}

In this section we perform the limit as $\eps \to 0$. We start by collecting some estimates coming from the results obtained in the previous section.

\begin{lemma}

{Given $\bbu_{\eps}^0 = \bbu^0 \in H^2(\Om)^{n+1}$} independent of $\eps$, there exists a constant $C>0$, independent of $\eps$, such that for any $i= 0, \dots, n$ it holds
    \begin{align}
    & \int_0^T \|\partial_t u_{\eps,i}\|_{{H^{1}(\Om)'}}^2 \; dt \le C,  \label{(a)} \\
    & \int_0^T \into |\nabla u_{\eps,i}|^2 \; \d xdt \le C, \label{a1} \\
    & {\|B_{\eps}(u_{\eps,i})\|_{L^2(0,T; L^2(\Omega))} \le C}, \label{b1} \\
     & \int_0^T \into \into K_{\eps}(x,y) |\nabla u_{\eps,i}(x)- \nabla u_{\eps,i}(y)|^2 \; \d x\d ydt \le C, \label{d1} \\
    &  \sum_{0 \le j \ne i \le n} \int_0^T \into |{J}_{\eps, ij}|^2 \; \d xdt \le C. \label{e1}
    \end{align}
\end{lemma}

\begin{proof}
{We first note that with the additional assumption of the initial datum in $H^2(\Omega)$ {it holds}
\[
\sum_{i=0}^n \into \into K_\eps(x,y) |\nabla u_{i}^{0}(x)- \nabla u_{i}^{0}(y)|^2 \; \d x \d y \le C,
\]
{where the constant $C>0$ depends on} $\|\bbu^0\|_{H^2(\Om)^{n+1}}$, but {not on} $\eps$, see  \cite[Theorem {1}]{BBM}. Thus, in this case all estimates in Lemma~\ref{thm:estim} hold with constants independent of $\eps$.}
Equations \eqref{b} and  \eqref{tau-uep} yield \eqref{d1}.
{Bound \eqref{b1} follows from \eqref{d1} and Lemma \ref{lemma.new.B.eps}.}
From Lemma \ref{lemma:J} and equation \eqref{e} we infer \eqref{e1}.
{Finally, \eqref{(a)} follows from \eqref{a1}, \eqref{d1}, \eqref{e1} and \eqref{eq:nonlocal}.
}
\end{proof}

 We are finally in a position to prove Theorem \ref{theorem-local}.

\begin{proof}[Proof of Theorem \ref{theorem-local}]
In view of \eqref{(a)} and \eqref{a1} we infer  \eqref{(i)}, while \eqref{ii)} follows from \eqref{b1} {and  Lemma \ref{lem:laplace}}.
It remains to show that such $\bbu$ actually satisfies the requirements in Definition \ref{def-local}. 

We first observe that, thanks to \eqref{(i)}, there holds
    \begin{equation}
        \label{u-1}
        u_{\eps, i} \to u_i \quad \text{strongly in } L^2(0,T; L^2(\Om)),
    \end{equation}
while \eqref{d1} and \cite[Lemma 4]{DRST} imply that
    \begin{equation}
    \label{u-2}
    u_{\eps,i} \to u_i \quad \text{strongly in } L^2(0,T; H^1(\Om)).
    \end{equation}
From Lemma \ref{lem:laplace} it follows that $\Delta u_{i}\in L^2(0,T; L^2(\Omega))$, therefore standard elliptic regularity results yield $u_i\in L^2(0,T; H^2(\Om))$. This implies that (1$_L$) and (2$_L$) are satisfied. Condition (3$_L$) {follows via \cite[Prop.~23.23]{Zei12} from the boundedness of
$\bbu_{\eps}$ and $\partial_t \bbu_{\eps}$ in $L^2(0,T; H^1(\Om))^{n+1}$ and $L^2(0,T; H^1(\Om)')^{n+1}$, respectively.}
In order to verify condition (4$_L$) we note the following. Let 
$\varphi_i \in L^2(0,T; H^2(\Om))$ of the form $\varphi(x,t)=\mathcal{X}(x)\mathcal{T}(t)$, with $\mathcal{T}\in L^2(0,T)$, $\mathcal{X}\in H^2(\Omega)$ arbitrary. Thanks to \eqref{(i)} there holds
    \[
    \int_0^T \langle \partial_t u_{\eps,i}, \varphi_i \rangle_{H^1(\Om)', H^1(\Om)} \; dt \to \int_0^T \langle \partial_t u_{i}, \varphi_i \rangle_{H^1(\Om)', H^1(\Om)} \; dt.
    \]
Moreover, using again \eqref{(i)} gives
    \[
    \int_0^T \into \sum_{0 \le j \ne i \le n} L_{ij} (u_{\eps,j} \nabla u_{\eps,i}- u_{\eps,i} \nabla u_{\eps, j}) \cdot \nabla \varphi_i \; \d xdt \to \int_0^T \into \sum_{0 \le j \ne i \le n} L_{ij} (u_{j} \nabla u_{i}- u_{i} \nabla u_{j}) \cdot \nabla \varphi_i \; \d xdt,
    \]
    as well as 
    \[
    \int_0^T\int_\Om u_{\eps, i} u_{\eps, j} \nabla \left(\sum_{k=0}^n 
    (a_{ik}-a_{jk}) u_{\eps,k}\right)\cdot \nabla \varphi_i\;\d x\d t \to \int_0^T\int_\Om u_{i} u_{j} \nabla \left(\sum_{k=0}^n (a_{ik}-a_{jk}) u_{k}\right)\cdot \nabla \varphi_i\;\d x\d t.
    \]
Concerning the last term in \eqref{eq:nonlocal}, thanks to Lemma \ref{lemma:J} we have
    \[
    \begin{split}
   & \int_0^T \into \sum_{0 \le j \ne i \le n} L_{ij} J_{\eps, ij} \cdot \nabla \varphi_i= -\int_0^T \into \divergenz (u_{\eps,i} u_{\eps,j} \nabla \varphi_i) (q_{\eps,i}(u_{\eps})- q_{\eps,j}(u_\eps)).
   \end{split}
    \]
Combining \eqref{u-1}, \eqref{u-2} and \eqref{ii)} { and by observing that the same argument of Lemma \ref{lemma:J} yields,
 up to the extraction of a subsequence,  the convergence}
    {\[
    u_{\eps, i} u_{\eps,j} \nabla (q_{\eps,i}(\bbu_\eps)- q_{\eps,j}(\bbu_\eps)) \to {J}_{ij} \quad \text{weakly in } L^2(0,T; L^2(\Om)^d)
    \]
    as $\eps \to 0$, with 
    \[
J_{ij}=u_{ i} u_{ j} \nabla(\Delta u_{i}- \Delta u_{j}),
    \]
in the weak sense defined in Lemma \ref{lemma:newJ} below, } 
it follows that
    \[
    \int_0^T \into \sum_{0 \le j \ne i \le n} L_{ij} J_{\eps, ij} \cdot \nabla \varphi_i \to \int_0^T \into \sum_{0 \le j \ne i \le n} L_{ij} J_{ij} \cdot \nabla \varphi_i.
    \]
Therefore the weak formulation in (4$_L$) holds for $\varphi_i \in L^2(0,T; H^2(\Om))$ of the form $\varphi(x,t)=\mathcal{X}(x)\mathcal{T}(t)$, with $\mathcal{T}\in L^2(0,T)$, $\mathcal{X}\in H^2(\Omega)$ arbitrary. However, standard density arguments allow us to extend the formulation first to $\varphi(x,t)=\mathcal{X}(x)\mathcal{T}(t)$, with $\mathcal{T}\in L^2(0,T)$, $\mathcal{X}\in H^1(\Omega)$ arbitrary, then to $\varphi\in L^2(0,T; H^1(\Omega))$ arbitrary.
This completes the proof of (4$_L$) and of the theorem.
\end{proof}
We conclude this section by proving an analogue result to Lemma \ref{lemma:J} concerning the convergence of the maps $J_{\eps,ij}$.
\begin{lemma}
\label{lemma:newJ}
For all $i \ne j= 0, \dots, n$ there exists ${J}_{ij} \in L^2(0,T; L^2(\Om)^d)$ such that $ J_{ij}= u_{i} u_{j} \nabla (\Delta u_i- \Delta u_j)$ holds in the weak sense, that is,
    \[
    \int_0^T \into {J}_{ij} \cdot {\eta} \; \d xdt= - \int_0^T \into (\Delta u_i- \Delta u_j) \divergenz(u_{i} u_{j} {\eta}) \; \d x dt,
    \]
for every ${\eta} \in L^2(0,T; H^1(\Om)^d) \cap L^\infty ((0,T) \times \Om; \R^d) $, with ${\eta} \cdot {n}= 0$ on $\rand$ and such that, up to the extraction of a subsequence,
    \[
    u_{\eps, i} u_{\eps,j} \nabla (q_{\eps,i}(\bbu_\eps)- q_{\eps,j}(\bbu_\eps)) \to {J}_{ij} \quad \text{weakly in } L^2(0,T; L^2(\Om)^d)
    \]
    as $\eps \to 0$.
\end{lemma}
\begin{proof} From \eqref{u-2} and the uniform $L^\infty$ bounds for $\bbu_{\eps}$ we deduce that $u_{\eps,i}u_{j,\eps}\to u_i u_j$ strongly in $L^2(0,T; H^1(\Omega))$ for every $i,j=0,\ldots,n$. Furthermore from \eqref{q}, \eqref{ii)} we get
\begin{align*}
q_{\eps,i}(\bbu_\eps)-q_{\eps,j}(\bbu_\eps)\rightharpoonup
-\sum_{k=0}^n(c_{ik}-c_{jk})\Delta u_k\quad\mbox{weakly in }L^2(0,T; L^2(\Omega)).
\end{align*}
From Lemma \ref{lemma:J}, estimate \eqref{e1} and the above convergence relations we deduce:
\begin{align*}
&\int_0^T \into {J}_{ij} \cdot {\eta} \; \d x \d t = \lim_{\eps\to 0}
\int_0^T \into {J}_{\eps,ij} \cdot {\eta} \; \d x \d t \\ 
&= - \lim_{\eps\to 0}\int_0^T \into (q_{\eps,i}(\bbu_\eps)- q_{\eps,j}(\bbu_\eps)) \divergenz(u_{\eps,i} u_{\eps,j} {\eta}) \; \d x \d t
=\int_0^T \into \sum_{k=0}^n(c_{ik}-c_{jk})\Delta u_k \divergenz(u_{i} u_{j} {\eta}) \; \d x \d t
\end{align*}
for every $\eta\in C^\infty(\Omega)$.
\end{proof}

\section*{Acknowledgements}

GM acknowledges the support of the Deutsche Forschungsgemeinschaft (DFG) via grant GZ: MA 10100/1-1 project number 496629752. The work of ED was supported by the Austrian Science Fund (FWF) through projects 10.55776/F65, 10.55776/Y1292, 10.55776/P35359, and 10.55776/F100800, as well as from the OeAD-WTZ project CZ04/2019 (M\v{S}MT\v{C}R 8J19AT013). JFP thanks the DFG for support via the Research Unit FOR 5387 POPULAR, Project No. 461909888. 
NZ thanks the University of Augsburg for financial support via the research program “Forschungspotentiale besser nutzen!”.
The authors are thankful to the anonymous Referee for their careful reading and insightful comments on the manuscript.

\section*{Conflict of interest}
The authors have no competing interests to declare that are relevant to the content of this article.

\section*{Data availability statement}
There is no additional data for this manuscript.

\end{document}